\RequirePackage{fix-cm}
\documentclass[smallextended]{svjour3}       
\smartqed  
\usepackage{graphicx}
\usepackage{amsmath}
\usepackage{mathptmx}      
\usepackage{enumerate}
\usepackage[misc]{ifsym}
\usepackage{subcaption}
\usepackage{bm}
\usepackage{amsfonts}
\usepackage{algorithm}
\usepackage{algpseudocode} 
\usepackage{epstopdf}
\usepackage{multirow}
\usepackage{siunitx}
\journalname{Mathematical Geoscience}

\usepackage[authoryear]{natbib}
\bibpunct{(}{)}{}{a}{}{;}

\begin{document}

\title{Smooth digital terrain modelling in irregular domain using finite element thin plate splines and adaptive refinement}

\titlerunning{Finite element thin plate spline for DTM}

\author{L. Fang}


\institute{Lishan Fang (\Letter) \at
              School of Mathematical Sciences, Huaqiao University, Quanzhou 362021, China \\
              Tel.: +86-18965196069\\
              \email{fanglishan@hqu.edu.cn} 
}

\date{Received: date / Accepted: date}

\maketitle

\begin{abstract}
Digital terrain models (DTMs) are created using elevation data collected in geological surveys using varied sampling techniques like airborne lidar and depth soundings. This often leads to large data sets with different distribution patterns, which may require smooth data approximations in irregular domains with complex boundaries. The thin plate spline (TPS) interpolates scattered data and produces visually pleasing surfaces, but it is too computationally expensive for large data sizes. The finite element thin plate spline (TPSFEM) possesses smoothing properties similar to those of the TPS and interpolates large data sets efficiently. This article investigates the performance of the TPSFEM and adaptive mesh refinement in irregular domains. Boundary conditions are critical for the accuracy of the solution in domains with arbitrary-shaped boundaries and are approximated using the TPS with a subset of sampled points. Numerical experiments are conducted on aerial, terrestrial and bathymetric surveys. It is shown that the TPSFEM works well in square and irregular domains for modelling terrain surfaces and adaptive refinement significantly improves its efficiency. A comparison of the TPSFEM, TPS and compactly supported radial basis functions indicates its competitiveness in terms of accuracy and costs.


\end{abstract}


\section{Introduction}\label{sec:intro}

\begin{figure} [h]
		\centering
		\begin{subfigure}[b]{0.44\textwidth}
                \centering
			\includegraphics[scale=0.28]{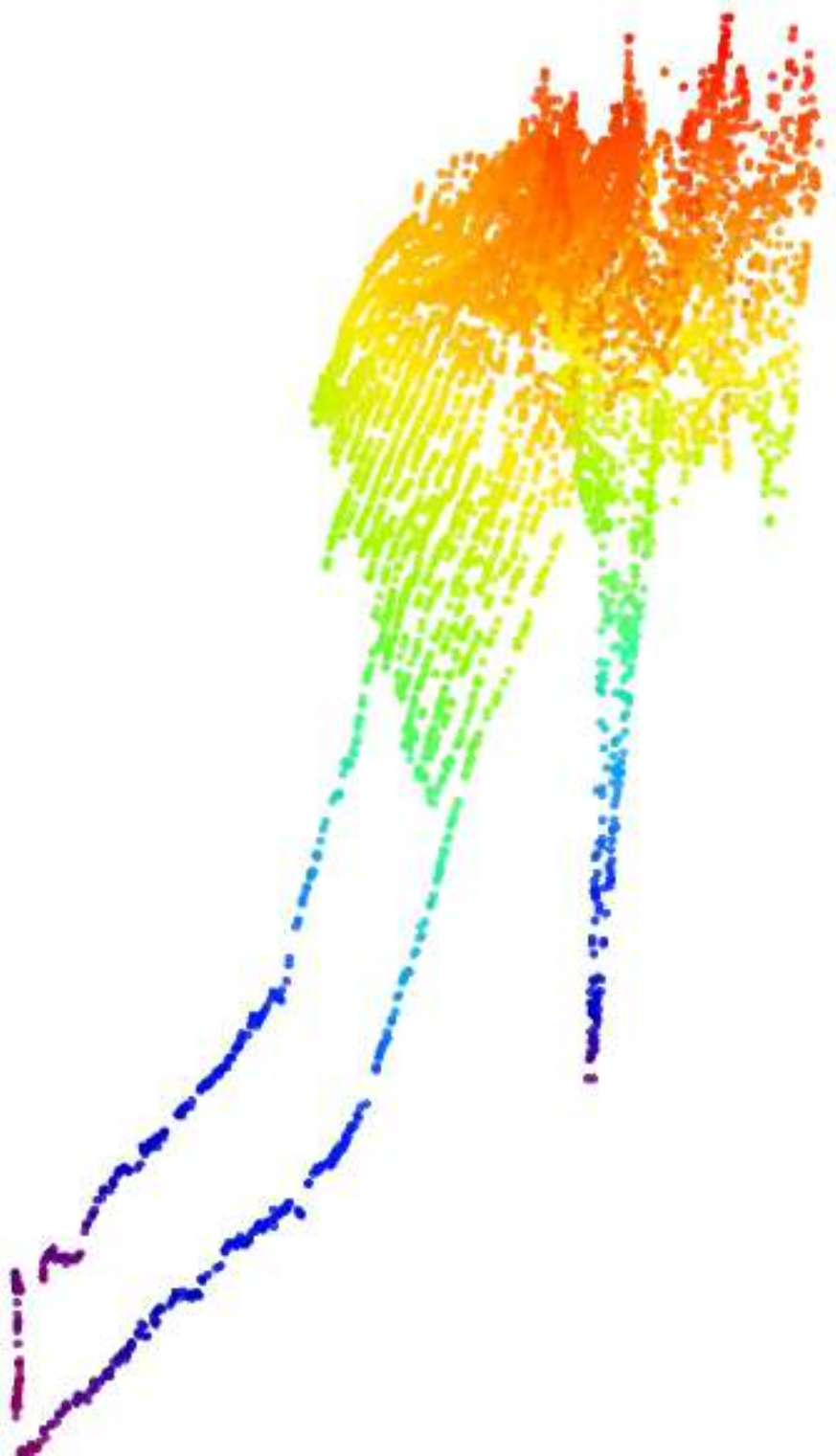}
			\caption{}
		\end{subfigure}
		\begin{subfigure}[b]{0.44\textwidth}
                \centering
			\includegraphics[scale=0.33]{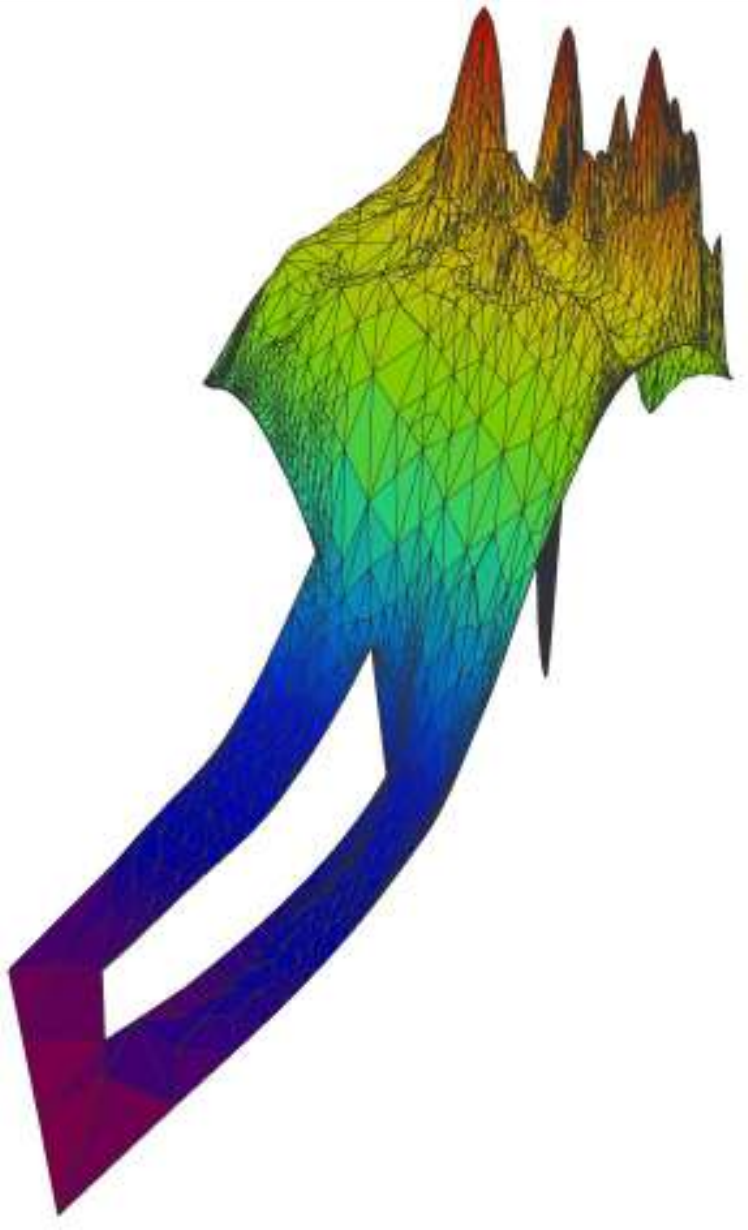}
			\caption{}
		\end{subfigure}
		\caption{(a) Observed data of a bathymetric survey; and (b) corresponding interpolant.}
		\label{fig:data_tpsfem}
	\end{figure}

Digital terrain models (DTMs) represent continuous terrain surfaces constructed using scattered data from geological surveys. They have been widely applied in various areas, including geological modelling, mining engineering and geomorphological mapping~\citep{galin2019review, mclaren2020visualisation}. The DTMs are normally built by interpolating observed data of terrain surfaces. An example is shown in Fig.~\ref{fig:data_tpsfem}, where irregularly distributed data are interpolated in an irregularly shaped domain with an interior hole. Popular interpolation techniques for the DTMs include kernel smoothing, wavelet-based methods, kringing, and smoothing splines, whose performance is influenced by the sizes and distribution patterns of data. Wavelet smoothing methods are most useful for data observed on a regular mesh in a rectangular domain and radial basis functions (RBFs) are applicable for arbitrarily scattered data~\citep{ramsay2002spline}. However, computational costs of RBFs like the thin plate spline (TPS) become impractical for large data sets. Methods like the random feature-based method~\citep{chen2020random} and triangulated irregular network~\citep{montealegre2015interpolation} were proposed to achieve higher efficiency, but they may not possess the flexibility and smoothing properties of the RBFs.


Many techniques were developed to improve the efficiency of the RBFs. Some use basis functions with compact support instead of global support for sparser systems of equations. \cite{wendland1995piecewise} presented several RBFs with compact support and \cite{deparis2014rescaled} introduced a scheme to adjust the radius of support for irregularly distributed data. Some others like the residual subsampling method produce smaller systems of equations~\citep{driscoll2007adaptive}. Another approach is to develop fast solvers for RBF interpolations. \cite{beatson1999fast} used fast matrix-vector multiplication and generalised minimal residual iteration with pre-conditioners to reduce the costs for calculating RBFs. 

Smooth interpolations based on second or higher-order finite elements (FEs) have also been developed for geological modelling. The finite element method (FEM) is often used to model complex models and works well in arbitrary topography. \cite{briggs1974machine} did some early work and proposed a discretised minimum curvature interpolation of the TPS to draw contour maps for the geophysical data. \cite{ramsay2002spline} introduced a discrete bivariate smoother based on quadratic FEs for smoothing over domains with complex irregular boundaries or holes. Recently, \cite{chen2018stochastic} approximated the TPS using nonconforming Morley FEs and provided stochastic convergence analysis. 

The finite element thin plate spline (TPSFEM) was developed by~\cite{roberts2003approximation} to interpolate and smooth large data sets of surface fitting and data mining applications. It is the first to use first-order FEs for a sparser system of linear equations that achieves higher efficiency~\citep{stals2006smoothing}. Adaptive mesh refinement and four error indicators were adapted from partial differential equations (PDEs) for the TPSFEM by~\cite{fang2024error}, which further enhanced its efficiency. Experiments showed that adaptive TPSFEM achieved lower root mean square error (RMSE) using less than $40\%$ of nodes and $50\%$ of system-solving time in square and L-shaped domains compared to uniform FE meshes.

\begin{figure} [h]
		\centering
		\begin{subfigure}[b]{0.45\textwidth}
                \centering
			\includegraphics[scale=0.8]{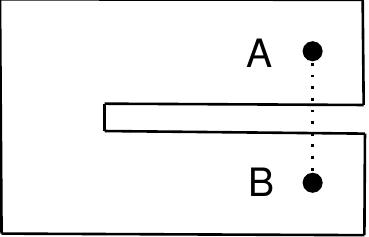}
			\caption{}
		\end{subfigure}
		\begin{subfigure}[b]{0.45\textwidth}
                \centering
			\includegraphics[scale=0.8]{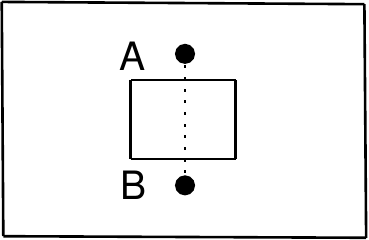}
			\caption{}
		\end{subfigure}
		\caption{(a) U-shaped domain; and (b) square domain with an interior hole.}
		\label{fig:domain_u_hole}
	\end{figure}

Irregular domains with complex boundaries are often required in applications like computer graphics, geometric modelling and geospatial modelling~\citep{brufau2002numerical, laga20183d}. They better reflect the geometry of the domain of interest and more accurately model surfaces compared to uniform square domains. However, irregular domains pose challenges to traditional RBFs, which rely on the Euclidean distance to measure the closeness between two points. While the RBFs work in domains with convex boundaries and no holes, they may smooth across concavities in the boundary or around the holes for irregular domains~\citep{ramsay2002spline}. An example is shown in Fig.~\ref{fig:domain_u_hole}, where lines connecting points A and B fall outside the domain. While these two points are close in Euclidean distance, they may not have a significant correlation. 


\begin{figure} [h]
		\centering
		\includegraphics[scale=0.33]{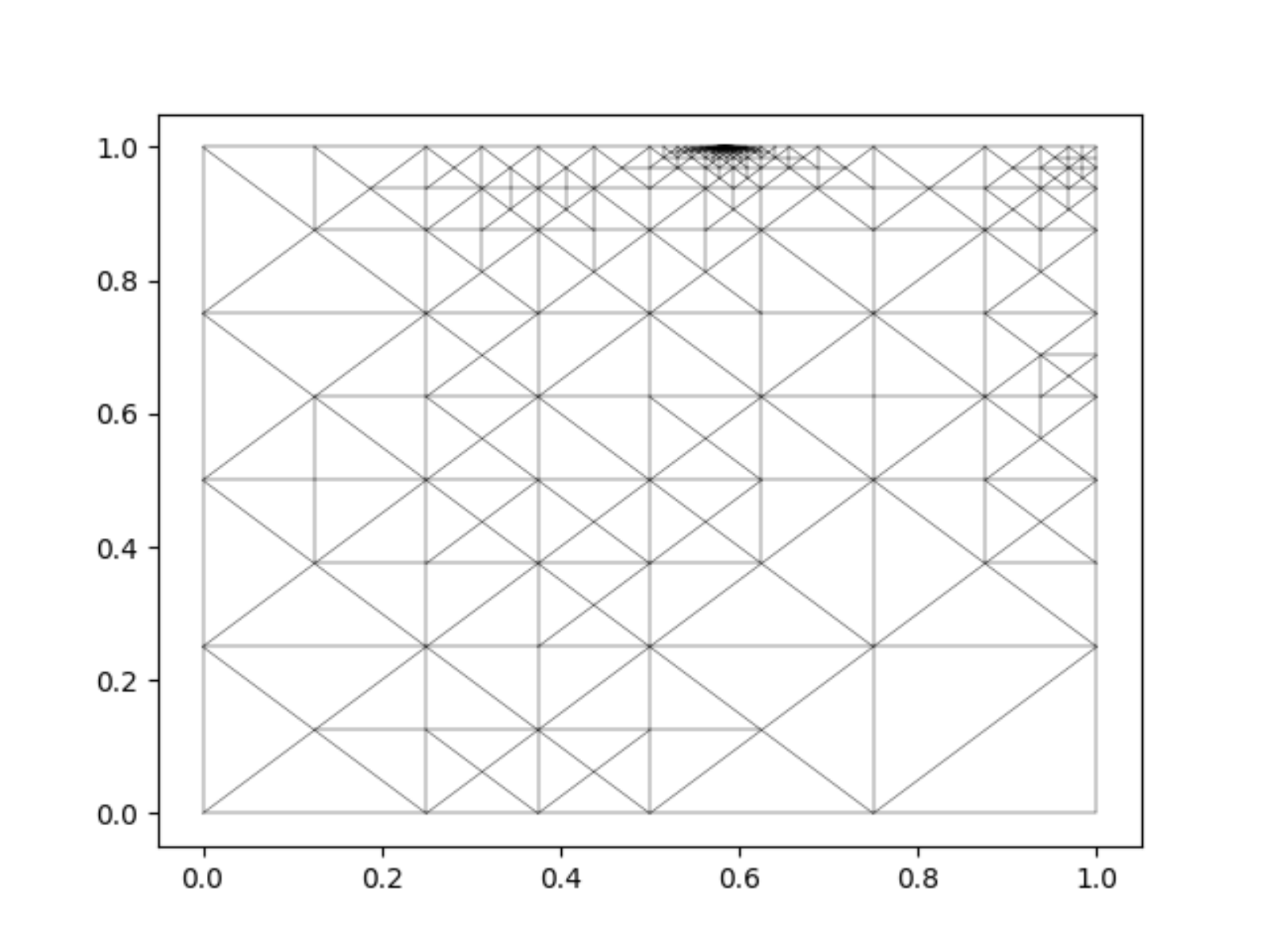}
		\caption{Adaptively refined mesh with over-refinement near $[0.6,1.0]$.}
		\label{fig:grid_overrefine}
	\end{figure}

This problem can be addressed using FEs to define the irregular topography. Ramsay's bivariate smoother interpolates on irregular domains but is expensive for large data sets~\citep{ramsay2002spline}. The TPSFEM interpolates large data sets efficiently and adaptive refinement markedly improves its efficiency in both square or L-shaped domains~\citep{fang2022adaptive,fang2024error}. However, the TPSFEM has two major problems with smoothing irregularly scattered data in irregularly shaped domains, which Ramsay stated is the most difficult practical scenario. Firstly, boundary conditions are critical for the accuracy of the TPSFEM and adaptive refinement but may be unknown for real-world data sets. Fang and Stals set Dirichlet boundaries as a constant number close to values of data, which will not work for irregular domains with varied data values close to boundaries. Accurate and efficient approximations of boundary conditions are required to build the TPSFEM in irregular domains. More information is given in Sect.~\ref{sec:boundary}. Secondly, over-refinement of FEs on oscillatory surfaces was found near the boundaries of FE meshes. An example is shown in Fig.~\ref{fig:grid_overrefine}, where FEs near $[0.6,1.0]$ are refined repeatedly while the other FEs remain coarse and the efficiency of the TPSFEM diminishes. This behaviour may become more severe in irregular domains with arbitrary-shaped boundaries and requires further investigation.



This article focuses on applications of the TPSFEM and adaptive mesh refinement for the DTMs in irregular domains. Dirichlet boundary conditions are approximated using the TPS with a small subset of sampled points. Over-refinement observed in previous studies becomes prevalent near complex boundaries in irregular domains and is alleviated by adjusting refinement routines. Numerical experiments are conducted on three real-world geological surveys to evaluate the performance of the TPSFEM and adaptive refinement in square and irregular domains. The analysis of the results focuses on the impacts of irregular domains and data distribution patterns on the TPSFEM and adaptive refinement. Lastly, a comparison of the TPSFEM to the TPS and two compactly supported RBFs is given to validate the TPSFEM's effectiveness.
 
The remainder of this article is organised as follows. Square and irregular domains are compared in Sect.~\ref{sec:domain}. The formulation and adaptive refinement of the TPSFEM is provided in Sect.~\ref{sec:tpsfem}. The Dirichlet boundary conditions approximated using the TPS are discussed in Sect.~\ref{sec:boundary}. The numerical experiments and data sets are described in Sect.~\ref{sec:experiment}. Results and discussions of the experiments are given in Sect.~\ref{sec:result}. Conclusions of findings are summarised in Sect.~\ref{sec:conclusions}.

\section{Square and irregular domains}\label{sec:domain}

\begin{figure} [h]
		\centering
		\begin{subfigure}[b]{0.49\textwidth}
                \centering
			\includegraphics[scale=0.4]{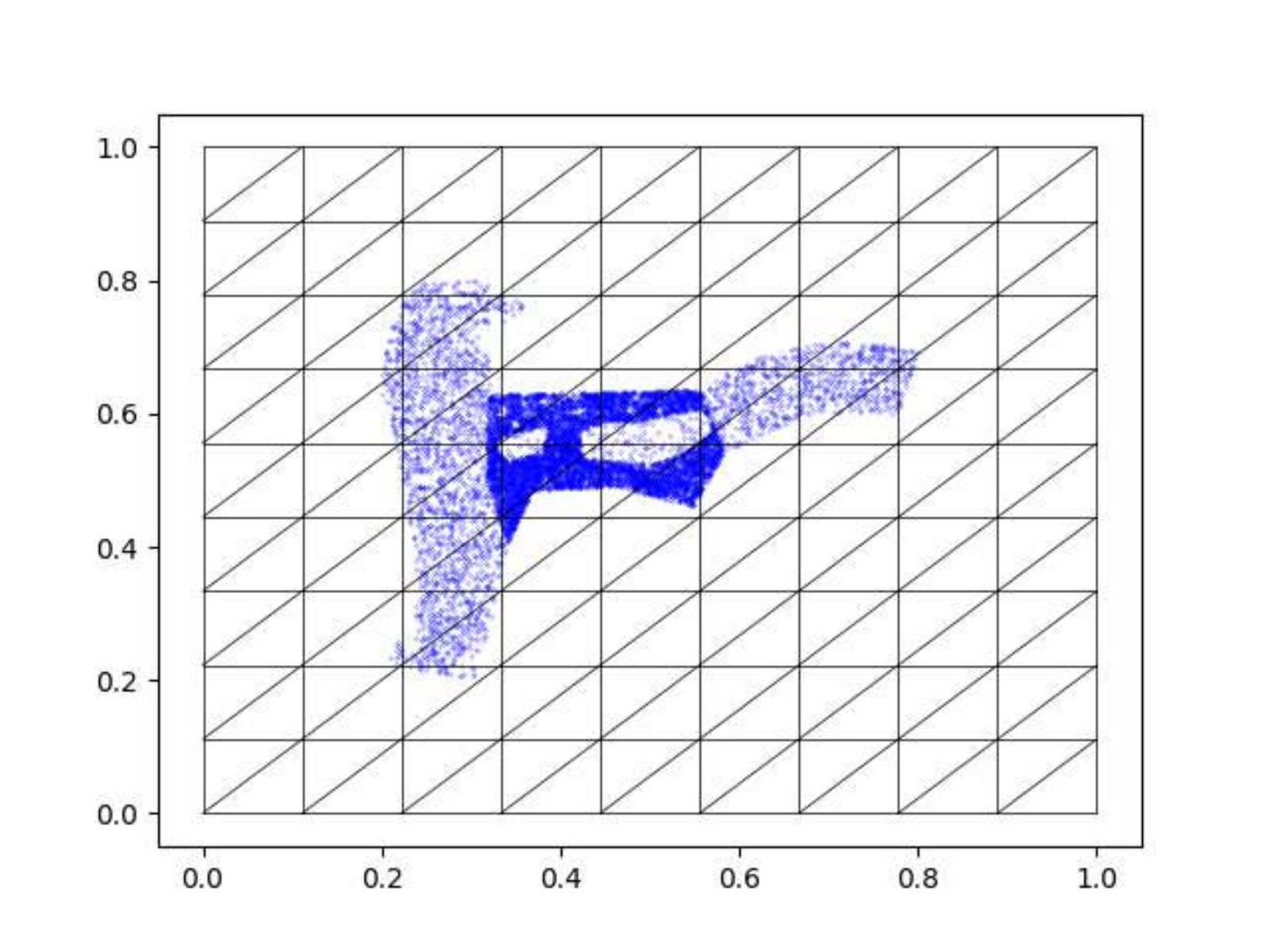}
			\caption{}
		\end{subfigure}
		\begin{subfigure}[b]{0.49\textwidth}
                \centering
			\includegraphics[scale=0.4]{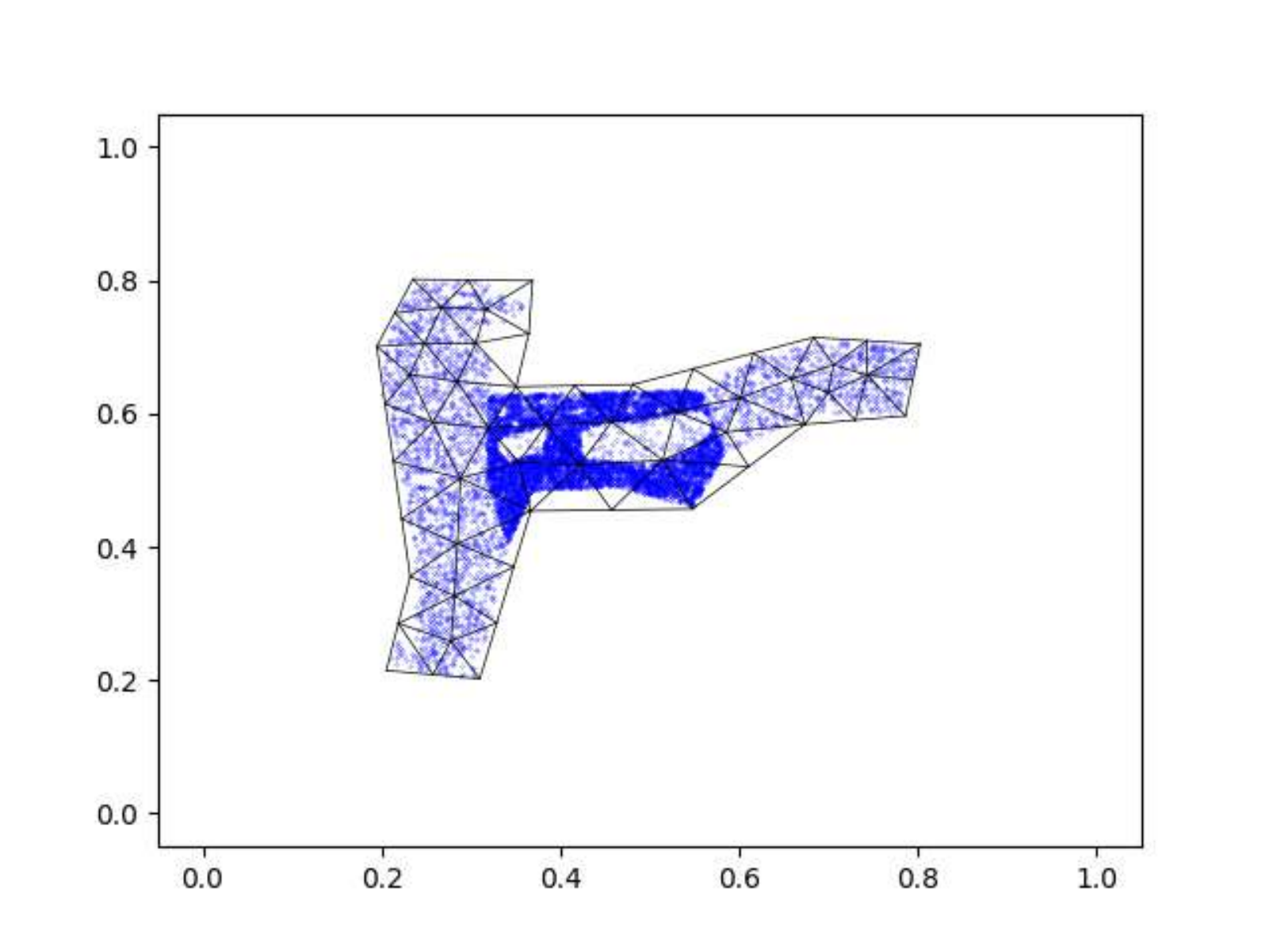}
			\caption{}
		\end{subfigure}
		\caption{FE meshes built in (a) square domain; and (b) irregular domain; with triangular elements. The data points are represented as blue dots.}
		\label{fig:mesh_square_irregular}
	\end{figure}

Recall that we are working with a data fitting technique based on FE discretisation. The first step of building the TPSFEM is constructing a FE mesh to contain the domain of the observed data. This article focuses on adaptive refinement over two-dimensional FE meshes with triangular elements and piecewise linear basis functions. The domain is partitioned into disjoint triangular elements and data points will be covered by at least one triangle. Examples of the observed data and corresponding FE meshes in square or irregular domains are shown in Fig.~\ref{fig:mesh_square_irregular}.

Square domains have been used to analyse properties of the TPSFEM and are applicable for data with various distribution patterns as shown in Fig.~\ref{fig:mesh_square_irregular}(a). However, they may lack accurate representations of the exterior shape of interpolants and cannot be used for modelling objects in complex geometric forms. Irregular domains fit the data closely and better reflect the domain of interest as shown in Fig.~\ref{fig:mesh_square_irregular}(b). The shape of irregular domains depends heavily on the configuration of data points, which often vary across applications. Aerial surveys collect data through flights that form frames and are less affected by terrain. In contrast, ground or bathymetric surveys may be limited by obstacles like the bank of a water body and sampled data have complex external geometric shapes. The construction of FE meshes in irregular domains is described in Appendix~\ref{app:domain}.

A potential benefit of irregular domains is the higher efficiency compared to square domains. Many elements in the FE mesh in Fig.~\ref{fig:mesh_square_irregular}(a) do not contain any data points and may contribute little to improving the accuracy of the solution. In contrast, most elements in Fig.~\ref{fig:mesh_square_irregular}(b) cover data points and the resulting interpolant may have higher efficiency. Note that this characteristic may overlap with adaptive mesh refinement and the impacts of their combinations are explored in Sect.~\ref{sec:result}.


\section{Finite element thin plate spline}\label{sec:tpsfem}

\subsection{Formulation}\label{sec:tpsfem_form}

Important concepts and notions of the TPS and TPSFEM are introduced in this section for further discussion. They have been explained in more detail in~\cite{roberts2003approximation}. The interpolation~$t$ of traditional TPS minimises functional
    \begin{equation} \label{eqn:tps_minimiser}
        J_{\alpha}(t) = \frac{1}{n}\sum^{n}_{i=1}\left(t(\bm{x}_{i})-y_{i}\right)^{2}+\alpha \int_{\Omega}\sum_{|\bm{v}|=2}\begin{pmatrix}2\\\bm{v}\end{pmatrix}\left(D^{\bm{v}}t(\bm{x})\right)^{2}\,d\bm{x}
    \end{equation}
over~$H^{2}(\Omega)$, where~$\bm{v}=(v_{1},\ldots,v_{d})$, $|v|=\sum_{j=1}^{d}v_{j}$ and~$v_{j}\in\mathbb{N}$ for~$j=1,\ldots,d$. Smoothing parameter~$\alpha$ balances the conformity to data and smoothness of~$t$. The construction of the TPS and an evaluation of its interpolant at a given point require $O(n^{3})$ and $O(n)$ arithmetic operations, respectively. As $n$ increases, the TPS becomes more accurate but computationally expensive.

The TPSFEM is defined on a FE mesh~$\mathcal{T}$ of dimension $d$, which covers the domain of interest~$\Omega$. Given observed data $\{(\bm{x}_{i}, y_{i}):i=1,2,\ldots,n\}$ of size~$n$, where~$\bm{x}_{i} \in \mathbb{R}^{d}$ is~$i$-th~$d$-dimensional predictor value and and~$y_{i} \in \mathbb{R}$ is~$i$-th response value. The TPSFEM~$s(\bm{x})=\bm{b}(\bm{x})^{T}\bm{c}$ is represented as a combination of~$m$ piecewise linear basis functions~$\bm{b}$, where~$\bm{b}=[b_{1},\ldots,b_{m}]^{T}$ and~$\bm{c}$ are corresponding coefficients. The basis functions $b_{k}$ corresponds to $k$-th node and $b_{k}(\bm{x})=0$ if $\bm{x}$ lies outside elements that contains the $k$-th node. 

    
    
 
The minimisation problem of the TPSFEM was adapted from Eq.~\eqref{eqn:tps_minimiser}, which requires continuous high-order derivatives to measure the smoothness of interpolants. While first-order derivatives $\nabla s$ can be computed using $s$ directly, the approximation of its second-order derivatives will be inaccurate~\citep{zienkiewicz2005finite}. \cite{roberts2003approximation} introduced auxiliary functions~$\bm{u}=[u_{1},\ldots,u_{d}]^{T}$ to represent the gradient of~$s$, where~$u_{j}$ approximate gradients of~$s$ in dimension~$j$ for~$1 \le j \le d$. Function $u_{j}$ is defined as $u_{j}(\bm{x})=\bm{b}(\bm{x})^{T}\bm{g}_{j}$ and~$\bm{g}_{j}$ are corresponding coefficients. 
 
 
The TPSFEM~$s$ minimises
	\begin{align} \label{eqn:minimiser}
		J(\bm{c},\bm{g}_{1},\ldots,\bm{g}_{d}) 
		& = \frac{1}{n}\sum^{n}_{i=1}\left(\bm{b}(\bm{x}_{i})^{T}\bm{c}-y_{i}\right)^{2}+\alpha \int_{\Omega} \sum_{j=1}^{d} \nabla \left(\bm{b}\left(\bm{x}\right)^{T}\bm{g}_{j}\right)\nabla \left(\bm{b}\left(\bm{x}\right)^{T}\bm{g}_{j}\right)\,d\bm{x} \nonumber \\
		& = \bm{c}^{T}A\bm{c}-2\bm{d}^{T}\bm{c}+\bm{y}^{T}\bm{y}/n+\alpha \sum_{j=1}^{d} \bm{g_{j}}^{T}L\bm{g_{j}},
	\end{align}
	subject to
	\begin{align} \label{eqn:constraint}
		L \bm{c}=\sum_{j=1}^{d}G_{j}\bm{g}_{j},
	\end{align}
 where~$\bm{d} = \frac{1}{n}\sum_{i=1}^{n}\bm{b}(\bm{x}_{i})y_{i}$, $\bm{y}=[y_{1},\ldots,y_{n}]^{T}$, $A = \frac{1}{n} \sum_{i=1}^{n}\bm{b}(\bm{x}_{i})\bm{b}(\bm{x}_{i})^{T}$. Matrix~$L$ is a discrete approximation to the negative Laplacian and $L_{pq}=\int_{\Omega}\nabla\bm{b}_{p}\nabla\bm{b}_{q}\,d\bm{x}$. Matrix~$G_{j}$ is a discrete approximation to the gradient operator in dimension~$j$ for~$1 \le j \le d$ and~$(G_{j})_{pq}=\int_{\Omega}\bm{b}_{p}\partial_{j}\bm{b}_{q}\,d\bm{x}$. The information of data is projected onto the FE mesh and stored in matrix~$A$ and vector~$\bm{d}$.
 Eq.~\eqref{eqn:constraint} ensures that $\nabla s$ and $\bm{u}$ are equivalent in a weak sense. Smoothing parameter~$\alpha$ is calculated using the generalised cross-validation (GCV)~\citep{stals2006smoothing}. 

The minimisation problem in Eq.~\eqref{eqn:minimiser} subject to Eq.~\eqref{eqn:constraint} is rewritten as a system of linear equations using Lagrange multipliers. For example, the system in a two-dimensional domain is written as
	\begin{equation*} \label{eqn:system}
		\begin{bmatrix} A & \bm{0} & \bm{0} & L \\ \bm{0} & \alpha L & \bm{0} & -G_{1}^{T} \\ \bm{0} & \bm{0} & \alpha L & -G_{2}^{T} \\ L & -G_{1} & -G_{2} & \bm{0} \end{bmatrix}\begin{bmatrix} \bm{c} \\ \bm{g}_{1} \\ \bm{g}_{2} \\ \bm{w} \end{bmatrix} = \begin{bmatrix} \bm{d} \\ \bm{0} \\ \bm{0} \\ \bm{0} \end{bmatrix}-\begin{bmatrix} \bm{h}_{1} \\ \bm{h}_{2} \\ \bm{h}_{3} \\ \bm{h}_{4} \end{bmatrix},
	\end{equation*}
where~$\bm{w}$ is a Lagrange multiplier and~$\bm{h}_{1}$, $\bm{h}_{2}$, $\bm{h}_{3}$, $\bm{h}_{4}$ store Dirichlet boundary information. The system's dimension is independent of the number of data points~$n$ instead of the number of basis~$m$, which improves its scalability. 

The error convergence of the TPSFEM smoother $s$ are dependent on $\alpha$, mesh size $h$ and maximum distance between data points~$d_{X}$~\citep{fang2024error}. While refinement of FEs reduces $h$, it cannot change $d_{\bm{X}}$ and the RMSE of $s$ may not be reduced. For example, adaptively refining a small region near boundaries may lead to small $h$ with relatively large $d_{\bm{X}}$, where the RMSE of $s$ is not reduced. This is undesired as it increases computational costs and the efficiency diminishes.


\subsection{Iterative adaptive refinement process}\label{sec:adaptive}

The efficiency of the TPSFEM may be improved by adaptive refinement that iteratively refines elements in sensitive regions like peaks and jumps marked by error indicators. The resulting adaptively refined meshes achieve higher precision with fewer nodes than uniform meshes, which is critical for smoothing large data sets. 

Four error indicators were adapted from FE approximations for PEDs and we concentrate on the auxiliary problem and recovery-based error indicators in this article~\citep{fang2024error}. The auxiliary problem error indicator approximates the error of the solution using the difference between interpolant $s$ and a locally more accurate approximation. The recovery-based error indicator post-processes gradients of $s$ to estimate errors. The former accesses neighbouring data points directly and the latter only uses values in FE meshes. Detailed descriptions of adaptive refinement and the error indicators are given in~\cite{fang2024error} and a brief version is provided in Appendix~\ref{app:adaptive}.

An iterative adaptive refinement process was adapted for the TPSFEM as described in Algorithm~\ref{alg:tpsfem_adaptive}, where $\mathcal{T}_{k}$, $\alpha_{k}$, $s_{k}$ represent the FE mesh, smoothing parameter and TPSFEM smoother in $k$-th iteration, respectively. It does not use a preset error tolerance value and employs an inner loop from steps 5 to 8 to refine FEs until the number of nodes doubles.


	\begin{algorithm}
		\caption{Iterative adaptive refinement process of TPSFEM}
		\label{alg:tpsfem_adaptive}
		\begin{algorithmic}[1]
            \Statex \textbf{Input:} Data and initial FEM grid~$\mathcal{T}_{0}$
            \Statex \textbf{Output:} TPSFEM smoother $s_k$
		\State Build initial mesh $\mathcal{T}_{0}$, calculate optimal~$\alpha_{0}$ and build TPSFEM~$s_{0}$ on $\mathcal{T}_{0}$; Calculate RMSE of~$s_{0}$.
		\State $k \gets 1$.
		\While{RMSE is higher than error tolerance or is not reduced less than $10\%$ for two iterations}
			\State calculate error indicators for edges in mesh~$\mathcal{T}_{k}$.
			\While{size of~$\mathcal{T}_{k}$ is smaller than twice of size of~$\mathcal{T}_{k-1}$}
				\State Iteratively mark and refine a subset of elements with high error indicator values.
				\State Calculate error indicator values for new edges in $\mathcal{T}_{k}$.
			\EndWhile
			\State Calculate optimal~$\alpha_{k}$ and build TPSFEM~$s_{k}$ on mesh~$\mathcal{T}_{k}$; Calculate RMSE of~$s_{k}$.
			\State $k \gets k+1$.
		\EndWhile
		\end{algorithmic}
	\end{algorithm}


\section{Approximations of Dirichlet boundary conditions}\label{sec:boundary}

Both Neumann and Dirichlet boundaries were investigated in~\cite{roberts2003approximation, stals2006smoothing} to analyse the TPSFEM. Dirichlet boundaries are chosen in this study since they are efficient and more stable for oscillatory surfaces near boundaries~\citep{fang2024error}. The Dirichlet boundary conditions of TPSFEM smoother $s$ are defined using~$\bm{c}$,~$\bm{g}_{1}, \bm{g}_{2}$ and~$\bm{w}$ values evaluated at boundary nodes. They are used to build $s$ and initialise values of new nodes on domain boundaries. Values of~$\bm{c}$,~$\bm{g}_{1}$ and $\bm{g}_{2}$ are calculated using function $f$ and its first-order derivatives directly. \cite{stals2006smoothing} showed that Lagrange multiplier $\bm{w}$ can be approximated using the second-order derivatives of $f$ if $\triangle f=0$ near boundaries.

Real-world data sets often do not contain boundary information and these values cannot be computed directly. \cite{fang2024error} replace~$\bm{c}$ values as the average of predictor values~$\bm{y}$ near boundaries and set~$\bm{g}_1=0$,~$\bm{g}_2=0$ and~$\bm{w}=0$ for Dirichlet boundaries. While it works for data sets with smooth surfaces and similar~$\bm{y}$ values near boundaries, it is not applicable in general situations. If boundary conditions are not sufficiently accurate, they may affect the accuracy of~$s$ and mislead the error indicators during adaptive refinement, especially when data points are close to the boundary.



\subsection{Approximations using TPS}\label{sec:boundary_tps}

Many RBFs like Gaussian, Multiquadrics, and TPS, are continuously differentiable and have been used to approximate the surface of scattered data and its derivatives~\citep{karageorghis2007matrix}. The TPS has smoothing properties similar to the TPSFEM and flexibility for small data sets (e.g. hundreds of data points). Thus, it is chosen to approximate Dirichlet boundary conditions for the TPSFEM. Despite that the TPS may become unstable near boundaries as argued by~\cite{fornberg2002observations}, it can greatly improve the accuracy of the TPSFEM $s$ near boundaries compared to~\cite{fang2024error}. The values of the TPSFEM interpolant are listed in Table~\ref{tab:tps_deri} with corresponding functions and kernels of the TPS. 
Given two predictor values~$\bm{x}_{i}$ and~$\bm{x}_{j}$, $\bm{x}_{i1}$ represents the first entry of~$i$-th predictor value.


	
\begin{table}
    \centering
    \caption{TPSFEM values, model problem functions and  TPS kernels.}
    \label{tab:tps_deri}
    \begin{tabular}{ccc}
		\hline\noalign{\smallskip}
		TPSFEM & Function & TPS Kernel \\
		\noalign{\smallskip}\hline\noalign{\smallskip}
		$\bm{c}$ & $f$ & $r^{2}\log(r)$ \\
		\hline\noalign{\smallskip}
		$\bm{g}_{1}$ & $\frac{\partial f}{\partial x_1}$ & $(2\log{r}+1)(\bm{x}_{i1}-\bm{x}_{j1})$ \\
		\hline\noalign{\smallskip}
		$\bm{g}_{2}$ & $\frac{\partial f}{\partial x_2}$ & $(2\log{r}+1)(\bm{x}_{i2}-\bm{x}_{j2})$ \\
		\noalign{\smallskip}\hline
		$-\frac{\bm{w}}{\alpha}$ & $\frac{\partial^2 f}{\partial x_1x_1}+\frac{\partial^2 f}{\partial x_2x_2}$ & $-\log{r}-4$ \\
		\noalign{\smallskip}\hline
    \end{tabular}
\end{table}

\begin{figure} [h]
		\centering
		\begin{subfigure}[b]{0.32\textwidth}
                \centering
			\includegraphics[scale=0.27]{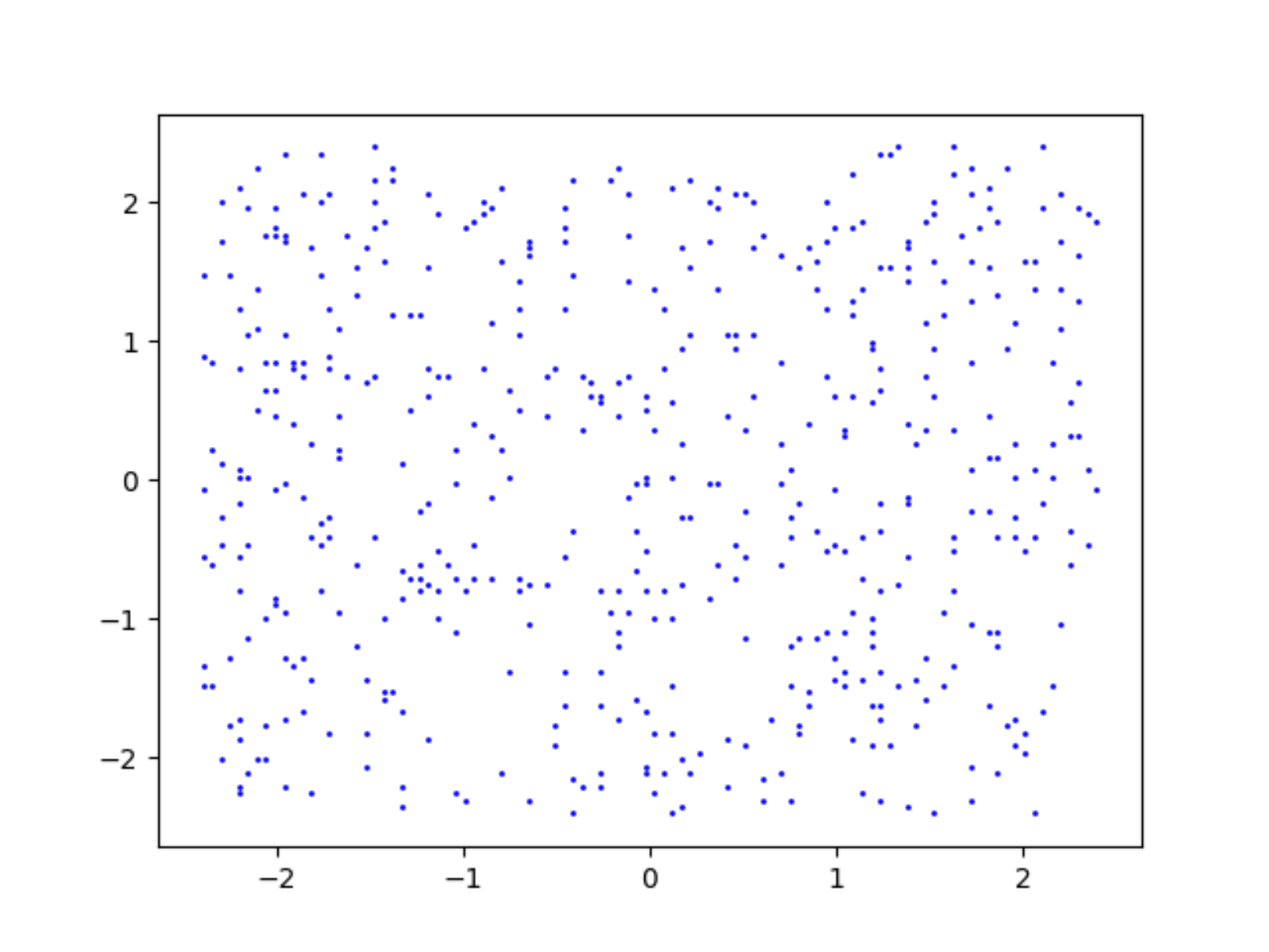}
			\caption{}
		\end{subfigure}
		\begin{subfigure}[b]{0.32\textwidth}
                \centering
			\includegraphics[scale=0.27]{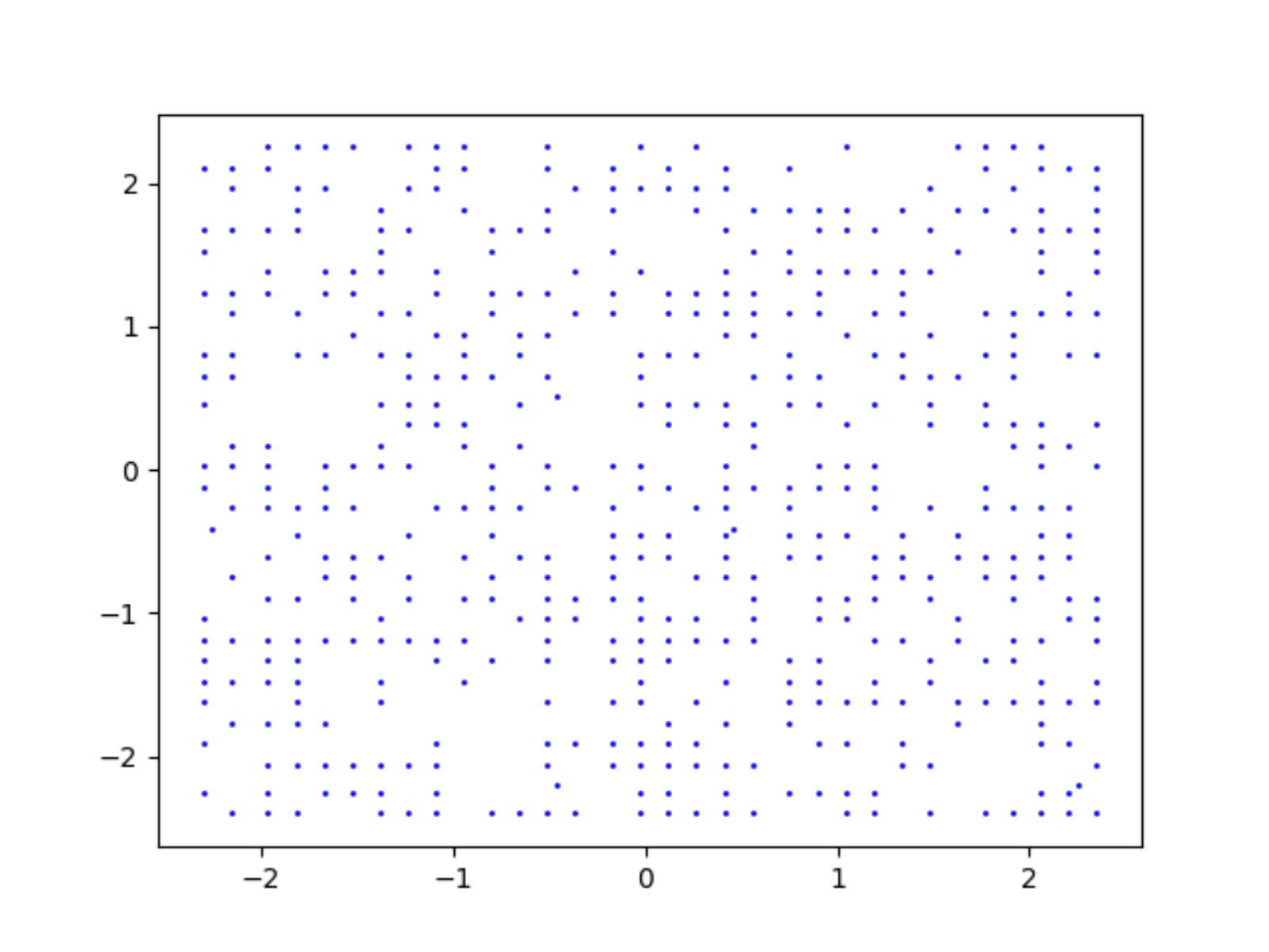}
			\caption{}
		\end{subfigure}
		\begin{subfigure}[b]{0.32\textwidth}
                \centering
			\includegraphics[scale=0.27]{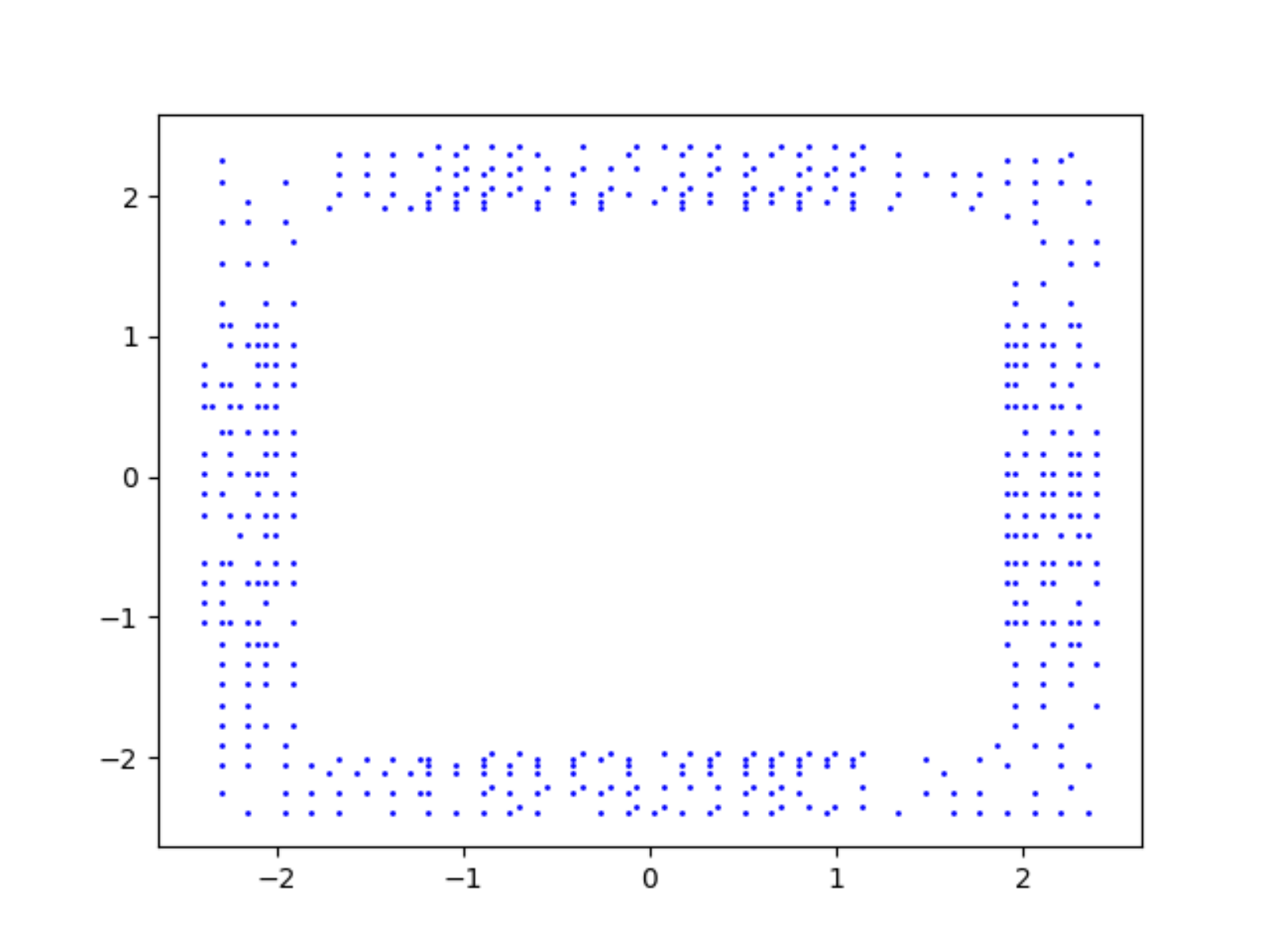}
			\caption{}
		\end{subfigure}
		\caption{Distribution of 500 data points within $[-2.4,2.4]^2$ region sampled (a) randomly; (b) using quadtree; and (c) using quadtree outside $[-1.9,1.9]^2$ region.}
		\label{fig:tps_data}
	\end{figure}

The accuracy of the TPS $t$ depends on a set of factors, including the choice of kernel, number $\hat{n}$ and distribution of sampled data. The number of sampled points $\hat{n}$ should be chosen as a trade-off between accuracy and costs and the sampling strategy may have large impacts especially when the sample size is small. A numerical experiment with noisy data is conducted to determine the optimal number and sampling strategy for $t$. Sampling data points randomly across the domain may lead to a large distance between neighbouring points as shown in Fig.~\ref{fig:tps_data}(a) and is not considered for the experiment. Alternatively, one may use quadtrees, which place rectangular meshes over the domain and randomly select data points from small rectangles. This helps to bound the maximum distance between neighbouring sampled points illustrated in Figs.~\ref{fig:tps_data}(b). Regions of sampling may also be restricted to areas near boundaries like in Fig.~\ref{fig:tps_data}(c) to improve its accuracy near boundaries. These two approaches are evaluated in the experiments described in Appendix~\ref{app:boundary_accuracy}.

The RMSE of approximations for $f$ and its first-order derivatives is reduced by increasing the number of basis and then it stabilises as the number reaches about 300 or 400. However, the second-order derivatives cannot be approximated accurately. While the boundary sampling performs relatively better than quadtree sampling in the experiment, the difference in RMSE gradually diminishes as $\hat{n}$ increases. Therefore, 300 data points are sampled to build the TPS using quadtrees for the numerical experiments in Sect.~\ref{sec:experiment}.

\subsection{Over-refinement near boundaries}\label{sec:overrefine}

Over-refinement is an undesired behaviour of adaptive refinement, where elements in certain regions are refined excessively and contribute little to improving the accuracy of the solution. This behaviour often occurs near abrupt changes like discontinuities, where finer elements may not increase the accuracy but still indicate high errors. \cite{fang2024error} showed that over-refinement also appears in oscillatory regions near boundaries and the resulting adaptive meshes underperform compared to uniform meshes. This weakens the stability of the TPSFEM and adaptive refinement.

\begin{figure} [h]
		\centering
		\begin{subfigure}[b]{0.49\textwidth}
                \centering
			\includegraphics[scale=0.4]{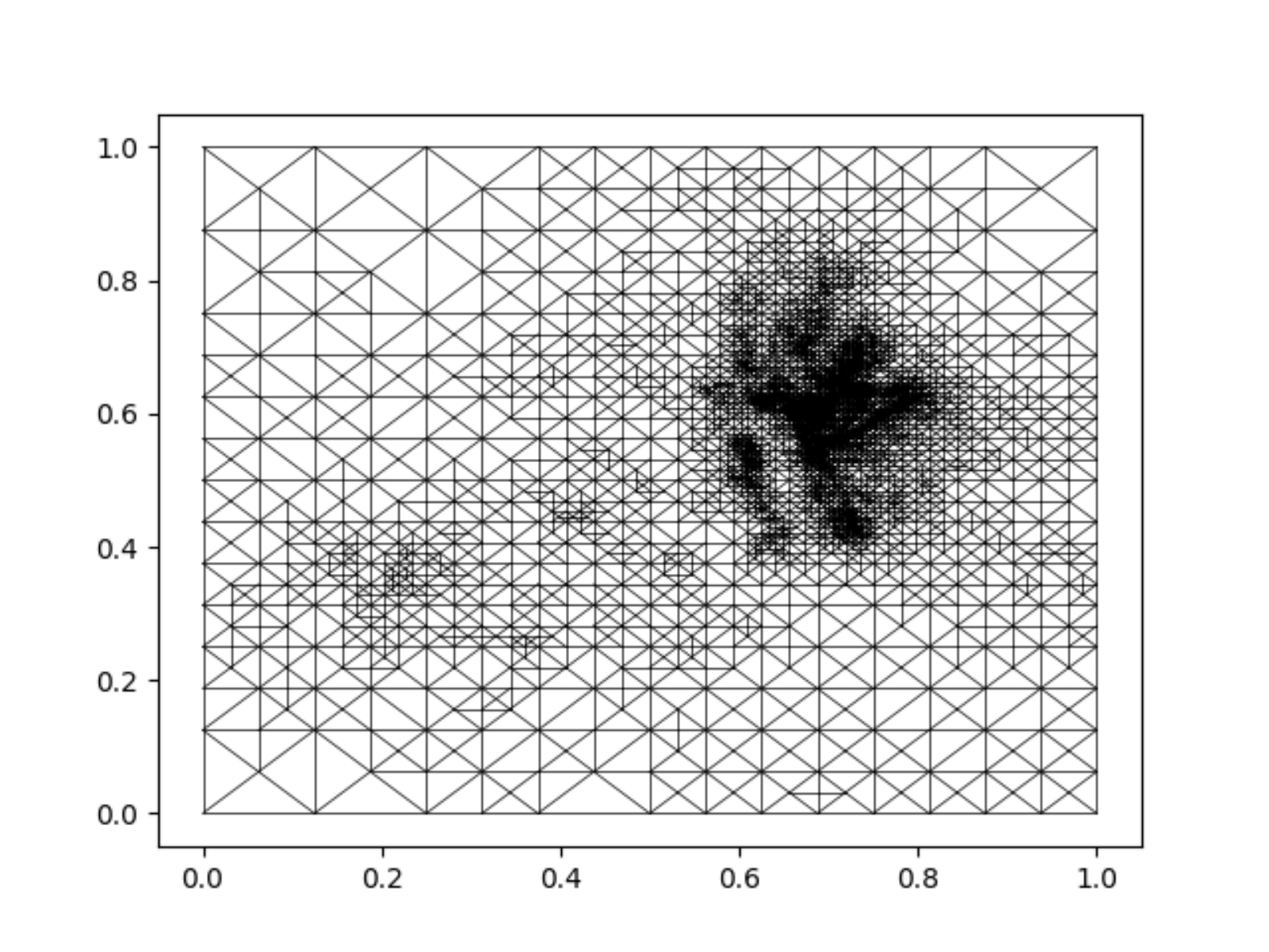}
			\caption{}
		\end{subfigure}
		\begin{subfigure}[b]{0.49\textwidth}
                \centering
			\includegraphics[scale=0.38]{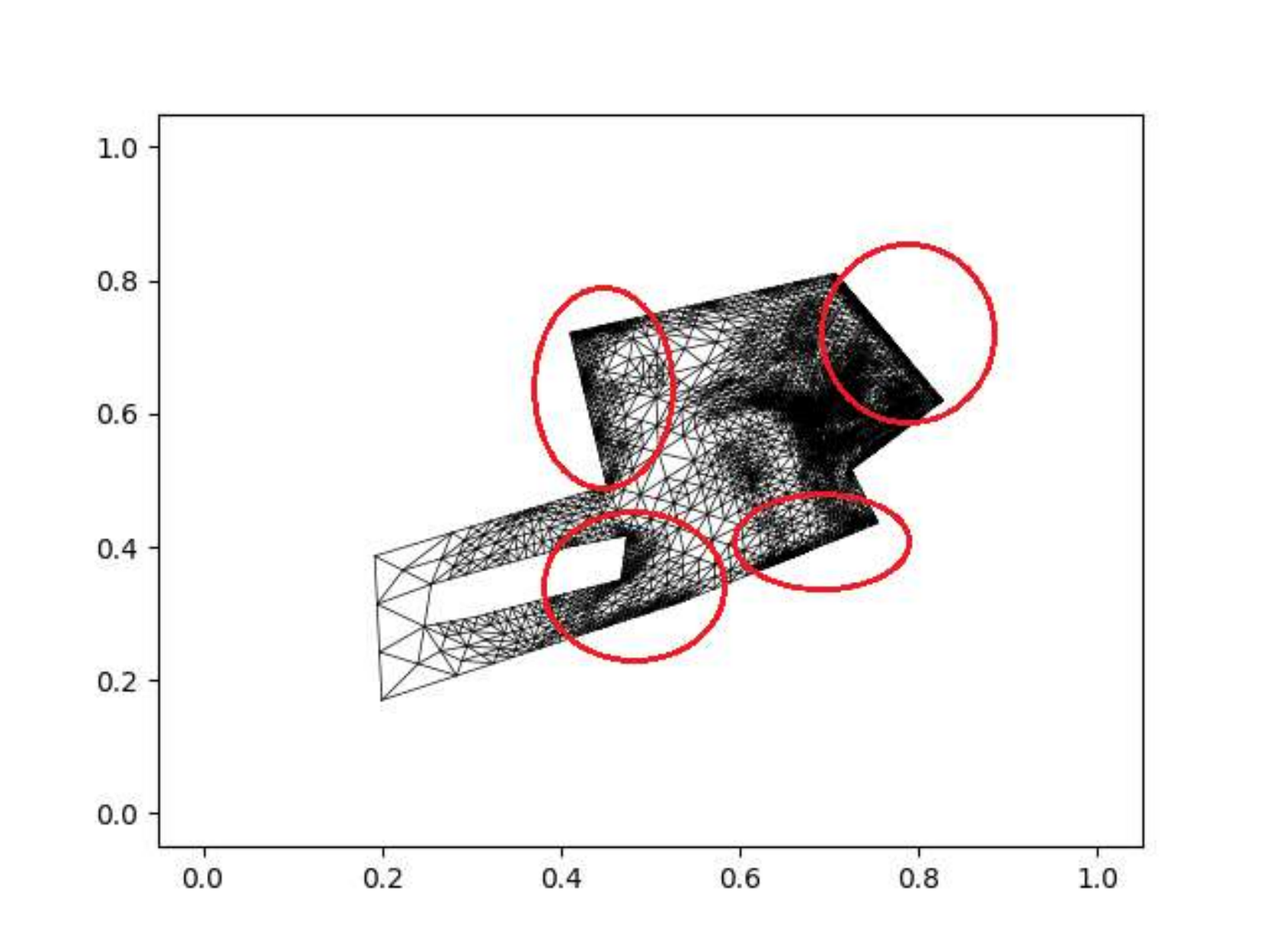}
			\caption{}
		\end{subfigure}
		\caption{FE meshes generated for River data using recovery-based error indicator in (a) square domain; (b) irregular domain. Four regions with over-refinement of FEs near boundaries are marked by red circles.}
		\label{fig:mesh_overrefine}
	\end{figure}

While the TPS approximations provide more accurate boundary conditions to build the TPSFEM smoother~$s$ like step 1 of Algorithm~\ref{alg:tpsfem_adaptive}, it may be insufficient for adaptive refinement. An example is shown in Fig.~\ref{fig:mesh_overrefine}, where two adaptive meshes are generated for the same data set in different domains and Dirichlet boundary conditions are approximated using the TPS. The finer elements of a square domain in Fig.~\ref{fig:mesh_overrefine}(a) concentrate on oscillatory regions at the centre of the domain. In contrast, the FE mesh in Fig.~\ref{fig:mesh_overrefine}(b) is constructed in an irregular domain and excessive refinement occurs near boundaries. They contribute little to improving the accuracy of the TPSFEM, which leads to lower efficiency.

Further investigations showed that most over-refinement occurs at step 6 of Algorithm~\ref{alg:tpsfem_adaptive}, where new nodes are initialised during adaptive refinement. For adaptive refinement of FE approximations for PDEs, if a new node is added to Dirichlet boundaries of FE meshes, its values are defined using known boundary functions like $f$. If a new interior node is added, its values are set as the average of neighbouring nodes. Since boundary conditions may be unknown for the TPSFEM, they must be approximated. Appendix~\ref{app:boundary_accuracy} shows that the RMSE stops decreasing as the number of basis reaches a certain number due to noise in data. The error indicators of the TPSFEM often rely on local information like a few data points or FEs and are sensitive to small changes. While the approximated boundary conditions are greatly improved compared to the constant values used in~\cite{fang2024error}, they may either be too expensive or insufficiently accurate for error indicators. This leads to oscillatory surfaces near boundaries, which may produce high error indicator values and excessive refinement of FEs near boundaries. 

An alternative for initialising new boundary nodes is proposed to avoid over-refinement during adaptive refinement. Instead of approximating new boundary node values during adaptive refinement, they are calculated as averaged values of their neighbouring nodes. Thus, surfaces near boundaries will not be altered and remain smooth, which does not indicate high errors. Note that the new boundary conditions are not more accurate and the accuracy of $s$ will not be improved in this process. Meanwhile, the TPS is still used to approximate Dirichlet boundary conditions to build $s$ for irregular domains.

These two approaches are tested and compared in the numerical experiments to evaluate their effectiveness. The original approach that uses TPS approximations to provide boundary conditions for $s$ and initialise new boundary node values during adaptive refinement is referred to as TPS approximation. The alternative that uses TPS approximations to build $s$ and initialise new boundary node values using the averaging is called nodal average.


 

\section{Numerical experiments}\label{sec:experiment}

The numerical experiments are conducted to evaluate the performance of the TPSFEM and adaptive refinement in both square and irregular domains. The three geological surveys for the experiments are from the United States Geological Survey\footnote{U.S. Geological Survey,https://www.usgs.gov/} and more details are given in Sect.~\ref{sec:data}. The settings of the experiments are provided in Sect.~\ref{sec:design}.

\subsection{Real-world data}\label{sec:data}
	
The first geological survey is an airborne magnetic survey collected from the Iron Mountain Menominee region of Michigan and Wisconsin~\citep{drenth2020airborne}. Precambrian rocks in the region are poorly mapped and the detailed high-resolution airborne magnetic survey is used to better understand its lithology and structure. It comprises 1,011,468 data points and variables including latitude, longitude and total magnetic field are selected for this study. The variations of magnetisation are related to differences in rock types.

The second survey contains ground observations collected by terrestrial laser scanners in Grapevine Canyon near Scotty's Castle, Death Valley National Park~\citep{morris2020geospatial}. It updated flood-inundation maps that were affected by a large flood in 2015 to current conditions. The observed data has been filtered of extraneous data such as vegetation, fences, power lines, and atmospheric interference for a continuous model. The resulting 122,104 data points, including latitude, longitude and elevation, are used to produce the DTM of the area.
	
The third survey is a nearshore bathymetric survey sampled at the mouth of the Unalakleet River in Alaska~\citep{snyder2019nearshore}. The survey comprises 895,748 data points consisting of latitude, longitude and height in meters regarding the reference ellipsoid. The data recorded depths from the seafloor to the echosounder and was collected using a small boat equipped with a single-beam sonar system.

	\begin{figure} [h]
		\centering
		\begin{subfigure}[b]{0.32\textwidth}
                \centering
			\includegraphics[scale=0.25]{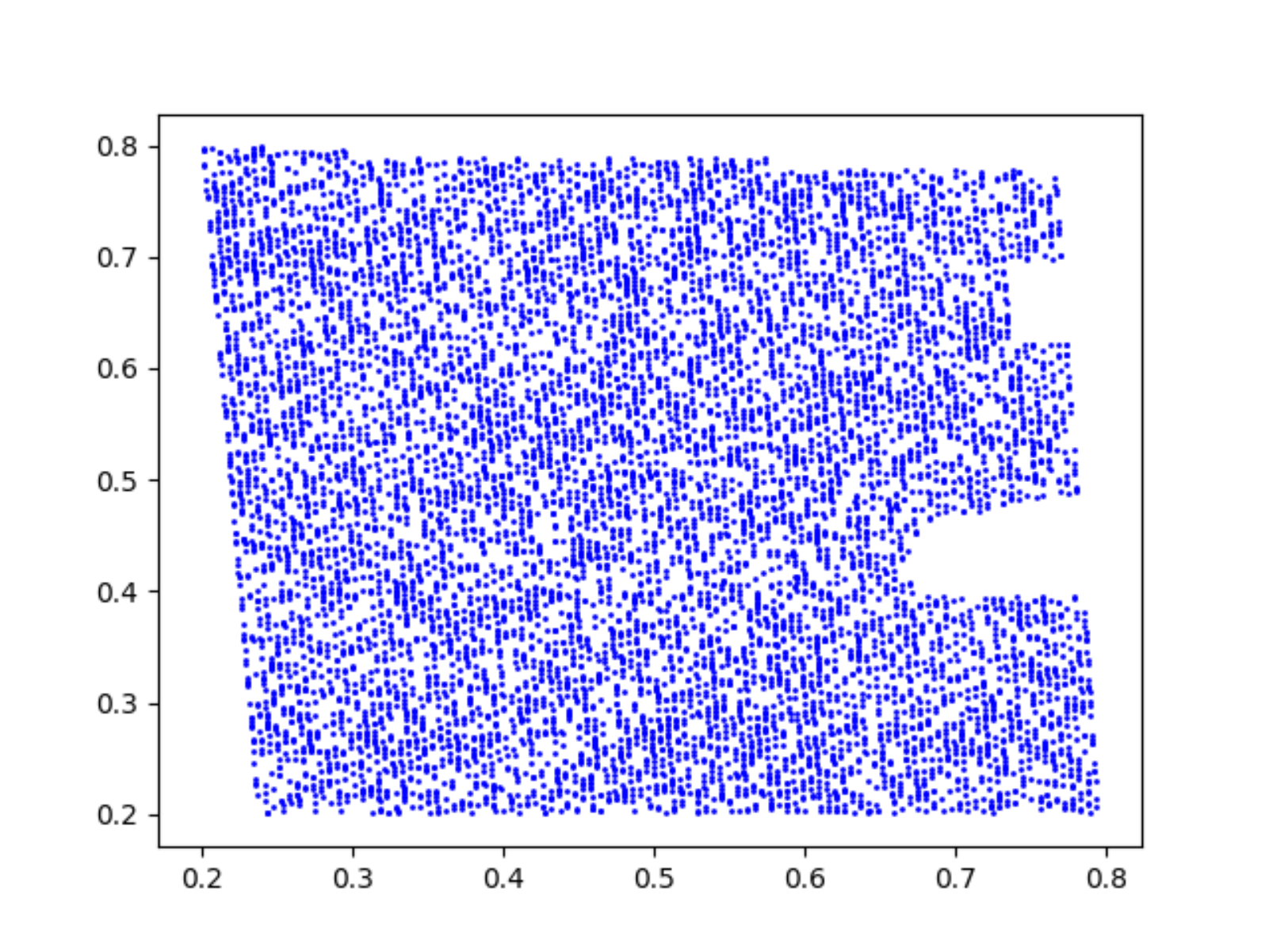}
			\caption{}
		\end{subfigure}
		\begin{subfigure}[b]{0.32\textwidth}
                \centering
			\includegraphics[scale=0.25]{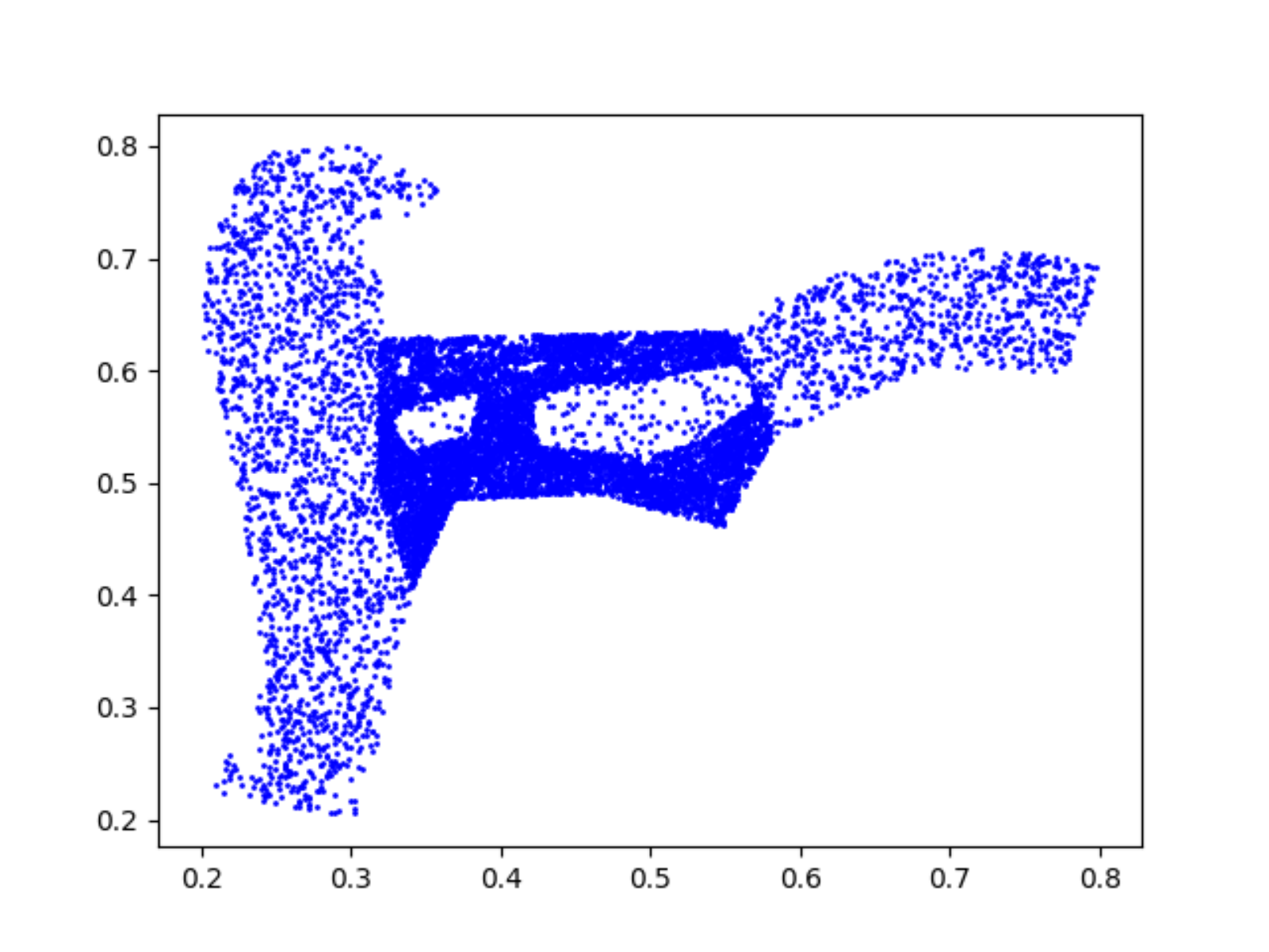}
			\caption{}
		\end{subfigure}
		\begin{subfigure}[b]{0.32\textwidth}
                \centering
			\includegraphics[scale=0.25]{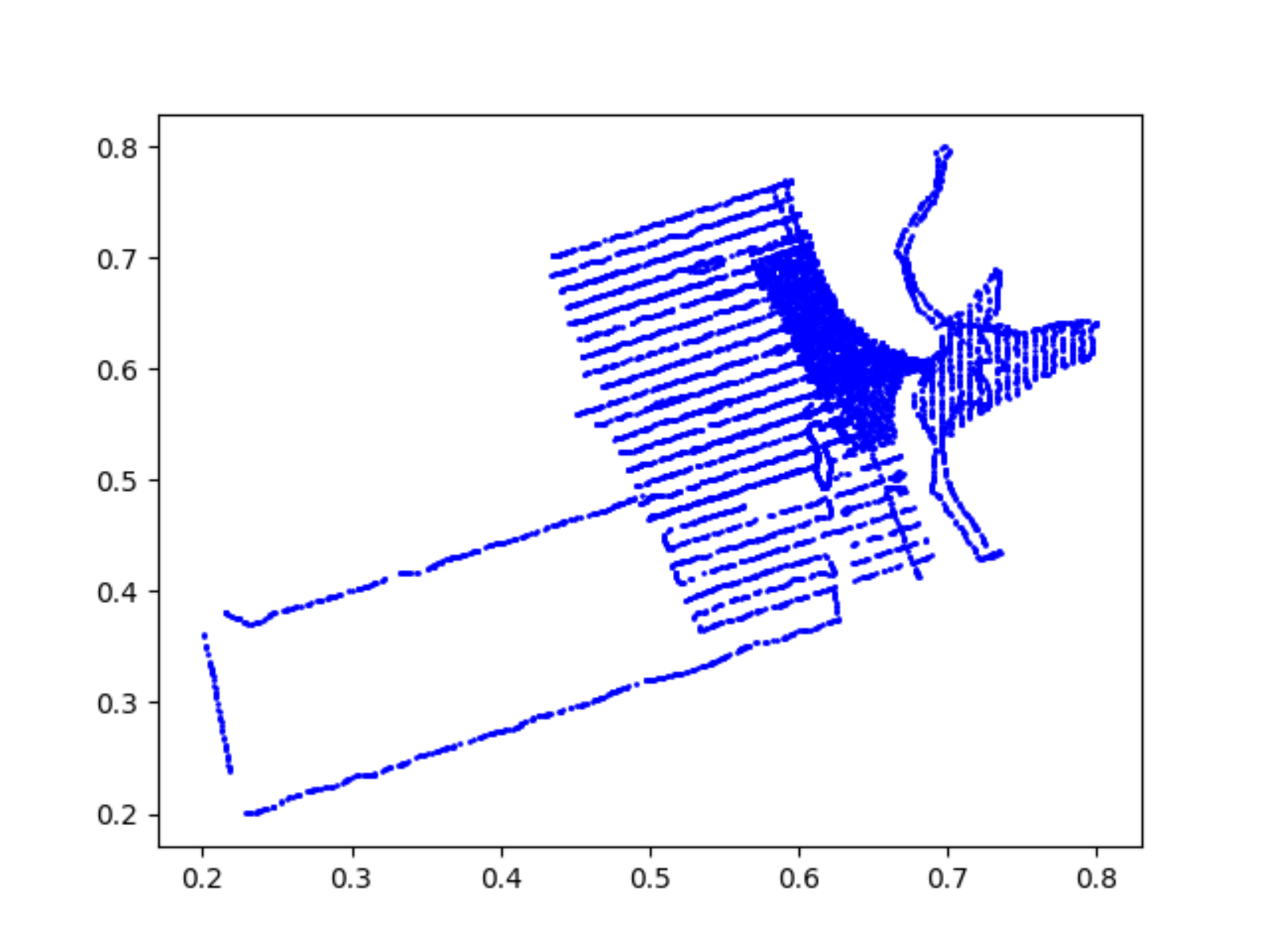}
			\caption{}
		\end{subfigure}
		\caption{Distributions of a subset of randomly selected $\bm{x}$ variables in (a) Mountain; (b) Canyon; and (c) River data sets.}
		\label{fig:data_dist}
	\end{figure}

	\begin{figure} [h]
		\centering
		\begin{subfigure}[b]{0.32\textwidth}
                \centering
			\includegraphics[scale=0.25]{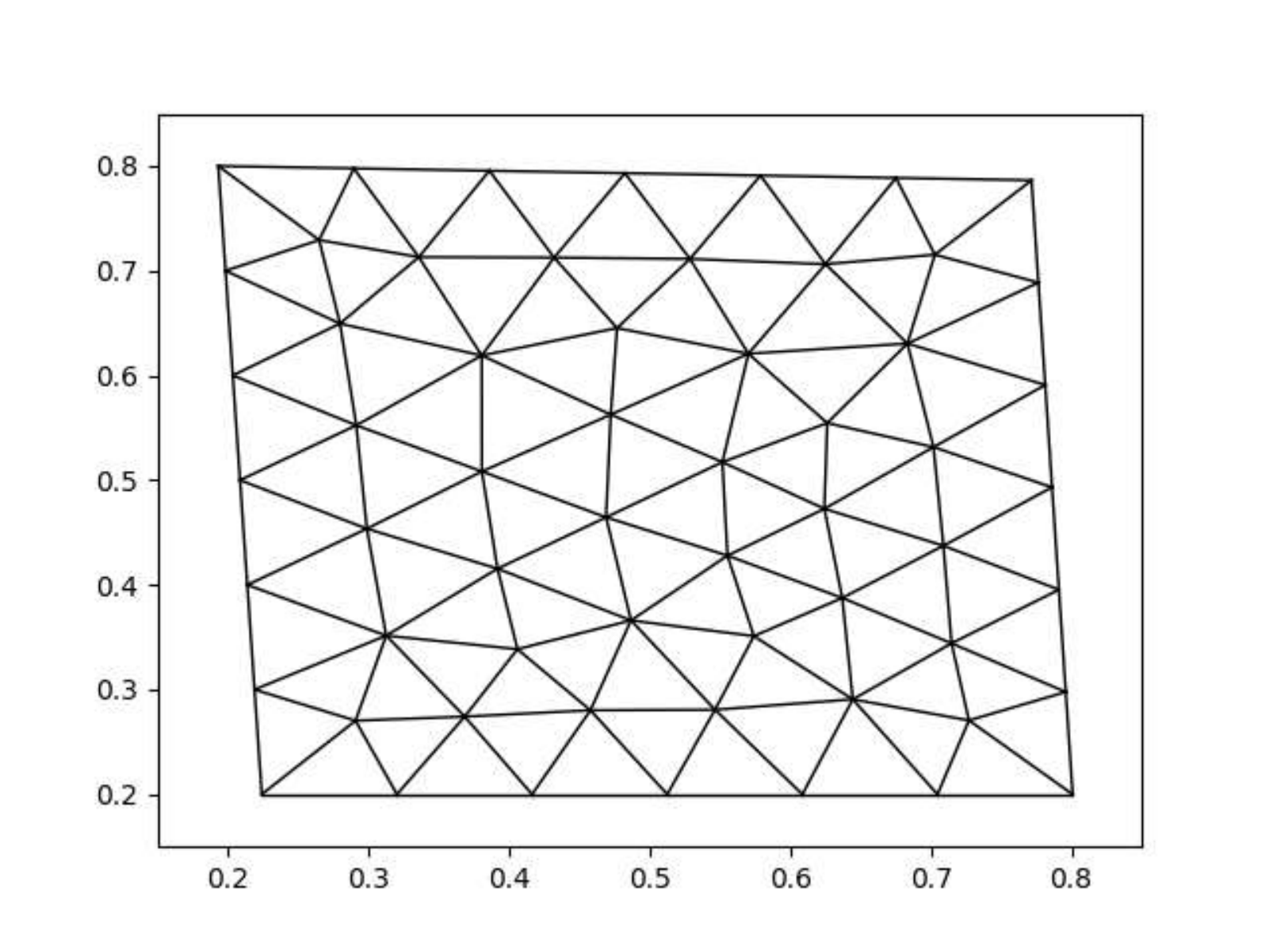}
			\caption{}
		\end{subfigure}
		\begin{subfigure}[b]{0.32\textwidth}
                \centering
			\includegraphics[scale=0.25]{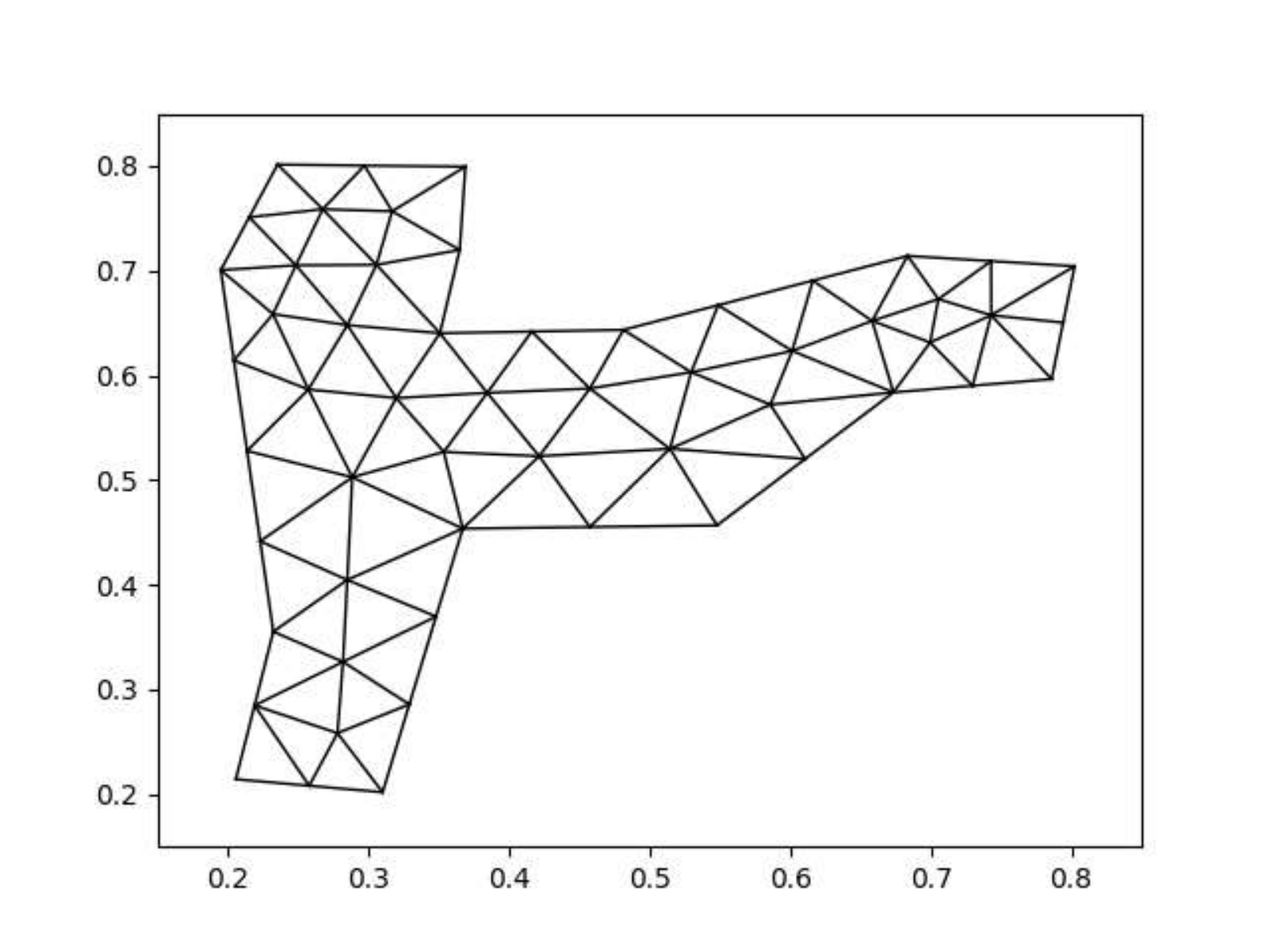}
			\caption{}
		\end{subfigure}
		\begin{subfigure}[b]{0.32\textwidth}
                \centering
			\includegraphics[scale=0.25]{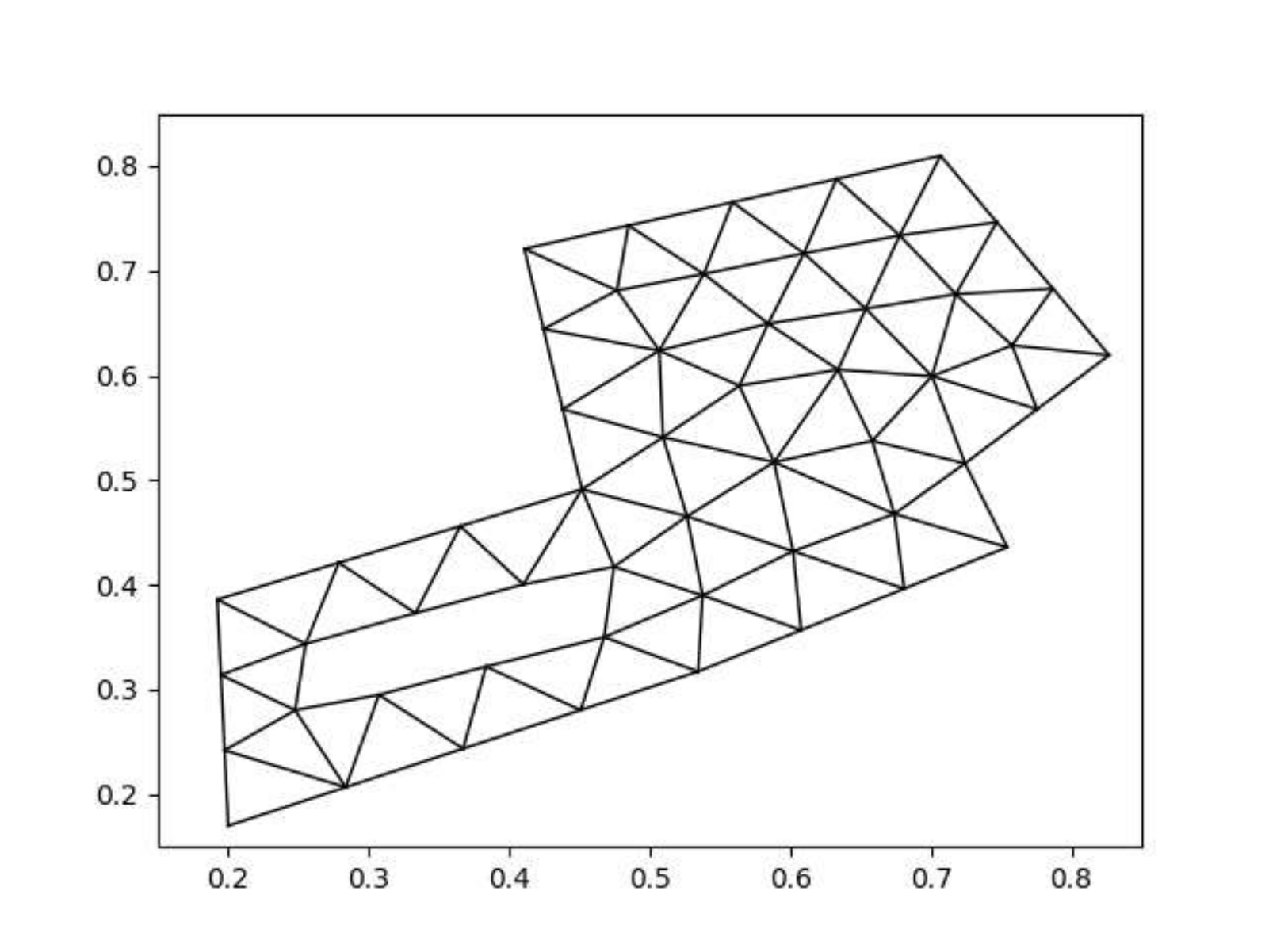}
			\caption{}
		\end{subfigure}
		\caption{Irregular FE meshes for (a) Mountain; (b) Canyon; and (c) River data sets.}
		\label{fig:data_domain}
	\end{figure}

The data from these three surveys will be referred to as the Mountain data, Canyon data and River data in the rest of the article. They are chosen due to their distinct data distribution patterns as illustrated in Fig.~\ref{fig:data_dist}. Corresponding initial irregular FE meshes are shown in Fig.~\ref{fig:data_domain}. The Mountain data was collected using an unmanned aerial vehicle in transects, which are relatively uniform across the domain as shown in Fig.~\ref{fig:data_dist}(a). In contrast, the Canyon data in Fig.~\ref{fig:data_dist}(b) are limited by terrains and equipment. Sampled points of the ground observations concentrate on the bottom and steep sides of the deep canyon that may be altered by the flood. More observations are located at the centre, which leads to uneven data densities. The River data in Fig.~\ref{fig:data_dist}(c) reflects routes of a boat that form frames with various sparseness. It also contains a region without any data point corresponding to theinterior holein the bottom left of the FE mesh in Fig.\ref{fig:data_domain}(c). More data points were sampled around the river mouth around point~$[0.65,0.6]$ to model complex river beds. 


\subsection{Experimental design}\label{sec:design}

Recall that we consider both square and irregular domains. All initial FE meshes in square domains contain 25 nodes and are refined using uniform or adaptive refinement for at most 10 and 8 iterations, respectively. The initial irregular FE meshes for the three data sets are generated using alphashape and pymesh as shown in Fig.~\ref{fig:data_domain} and contain 59, 56 and 53 nodes, respectively. They are refined using uniform or adaptive refinement for 8 or 7 iterations, respectively. The final square and irregular domains have similar numbers of nodes, which allows us to compare their performance. The predictor values~$\bm{x}$ of the three data sets are scaled and fitted inside the~$[0.2,0.8]^{2}$ region of a square~$[0,1]^{2}$ domain. The response values~$\bm{y}$ are sampled on different scales and normalised in the range of~$[0,1]$ for comparison. 

The efficiency of $s$ is measured using the RMSE versus the number of nodes in meshes, where RMSE~$ = \big(\frac{1}{n}\sum_{i=1}^{n}\big(s(\bm{x}_{i})-y_{i}\big)^{2}\big)^{\frac{1}{2}}$. Additionally, the ratio of interior nodes within the radius of 0.005 of boundary nodes indicates the degree of over-refinement near boundaries. This radius is chosen such that the ratio only includes nodes generated during over-refinement near boundaries. Note that the lengths of more than 50\% of edges in the final adaptive FE meshes are shorter than 0.005.
	
Their performance is compared through convergence plots and tables containing statistics of final meshes generated using different approaches. The plots display the convergence of RMSE in logarithmic scales for the three data sets using both uniform or adaptive refinement with the auxiliary problem and recovery-based error indicators described in Appendix~\ref{app:adaptive}. Their performance is evaluated in square and irregular domains with the TPS approximation or nodal average. Additional metrics of the final meshes are provided in tables, including the number of nodes, time for solving Eq.~\eqref{eqn:minimiser} of the TPSFEM, RMSE and maximum absolute residual errors (MAX). The computational costs, including time for building and solving systems, error indicators and GCV, of the TPSFEM are affected by factors like the number of nodes in FE meshes, number and distribution patterns of data points. A detailed analysis is too complicated as discussed in \cite{fang2024error} and not the focus of this article. Hence, we only include system-solving time, which reflects the difference in their efficiency.

The experiments are conducted using a TPSFEM program\footnote{TPSFEM program, https://bitbucket.org/fanglishan/tpsfem-program/} implemented in Python 3.8. The systems are solved by a sparse solver (PyPardiso\footnote{PyPardiso, https://pypi.org/project/pypardiso/}) and the time is measured in seconds. All experiments are conducted on a personal computer with Microsoft Windows 10 operating system and an Intel$^\circledR$ Core\textsuperscript{TM} i9-10900 CPU @ 2.80 GHz configuration with 64.0 GB of RAM.

\section{Results and discussions}\label{sec:result}

The numerical results and discussions for the Mountain, Canyon and River data are provided in Sects.~\ref{sec:moutain}, \ref{sec:canyon}, \ref{sec:river}, respectively. A comparison of the TPSFEM to the TPS and two compactly-supported RBFs (CSRBFs) is given in Sect.~\ref{sec:csrbf}.

\subsection{Mountain data}\label{sec:moutain}

\begin{figure} [h]
		\centering
		\includegraphics[scale=0.35]{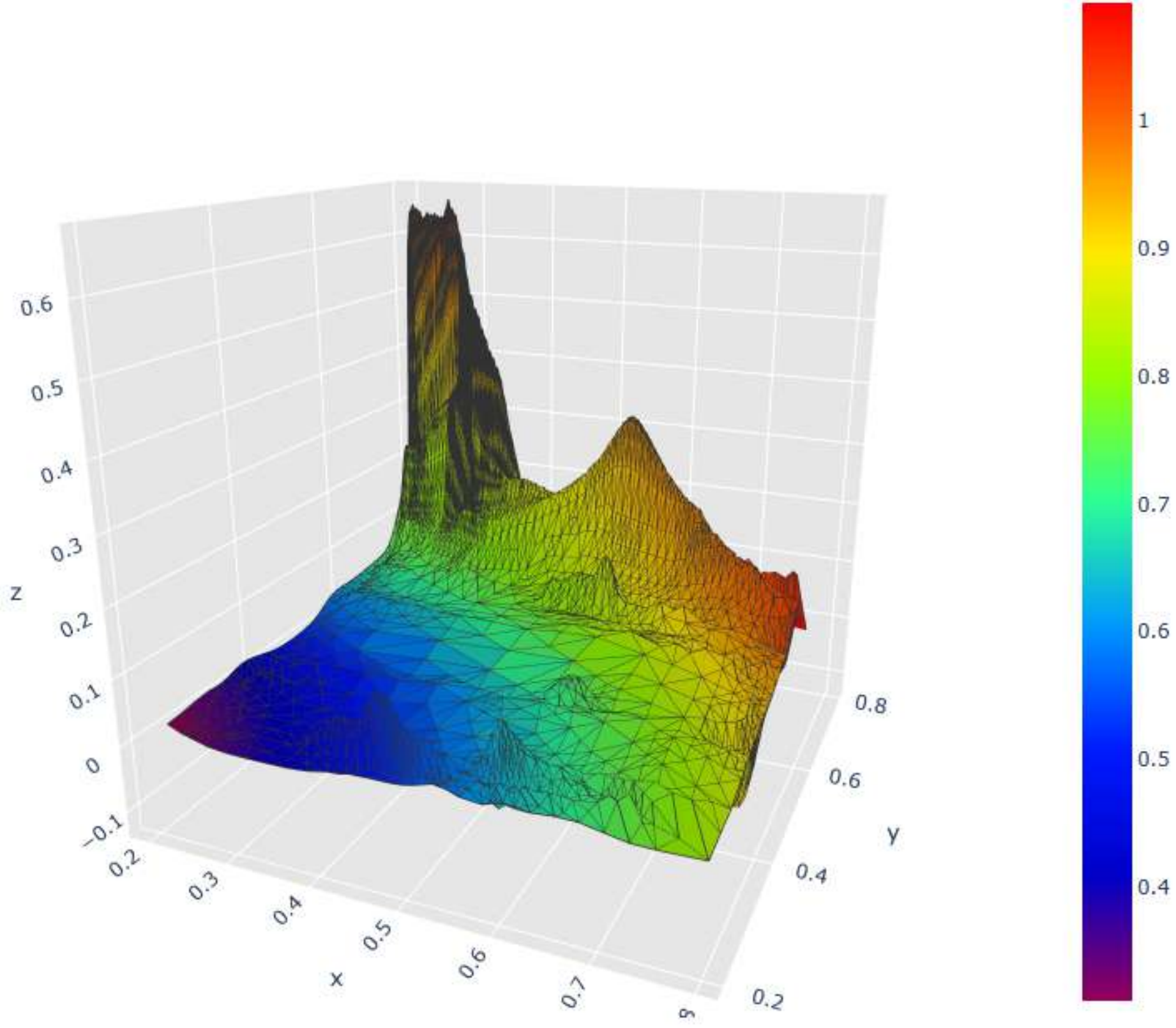}
		\caption{Triangulated surface of TPSFEM for Mountain data in an irregular domain.}
		\label{fig:trisurf_mountain}
	\end{figure}

A triangulated surface plot for the Mountain data generated using the recovery-based error indicator in an irregular domain is shown in Fig.~\ref{fig:trisurf_mountain}. The Mountain data consists of several oscillatory surfaces at the upper left corner and relatively smooth surfaces in the rest of the domain. High peaks near $[0.25,0.8]$ show high magnetic field readings and suggest potential mineral resources. Adaptive refinement accurately detects oscillatory regions and iteratively refines them to achieve the required accuracy. The resulting mesh contains markedly finer elements over the peaks and coarser elements elsewhere, which effectively reduces computational costs and memory requirements of the TPSFEM as discussed below.

\begin{figure} [h]
		\centering
		\includegraphics[scale=0.26]{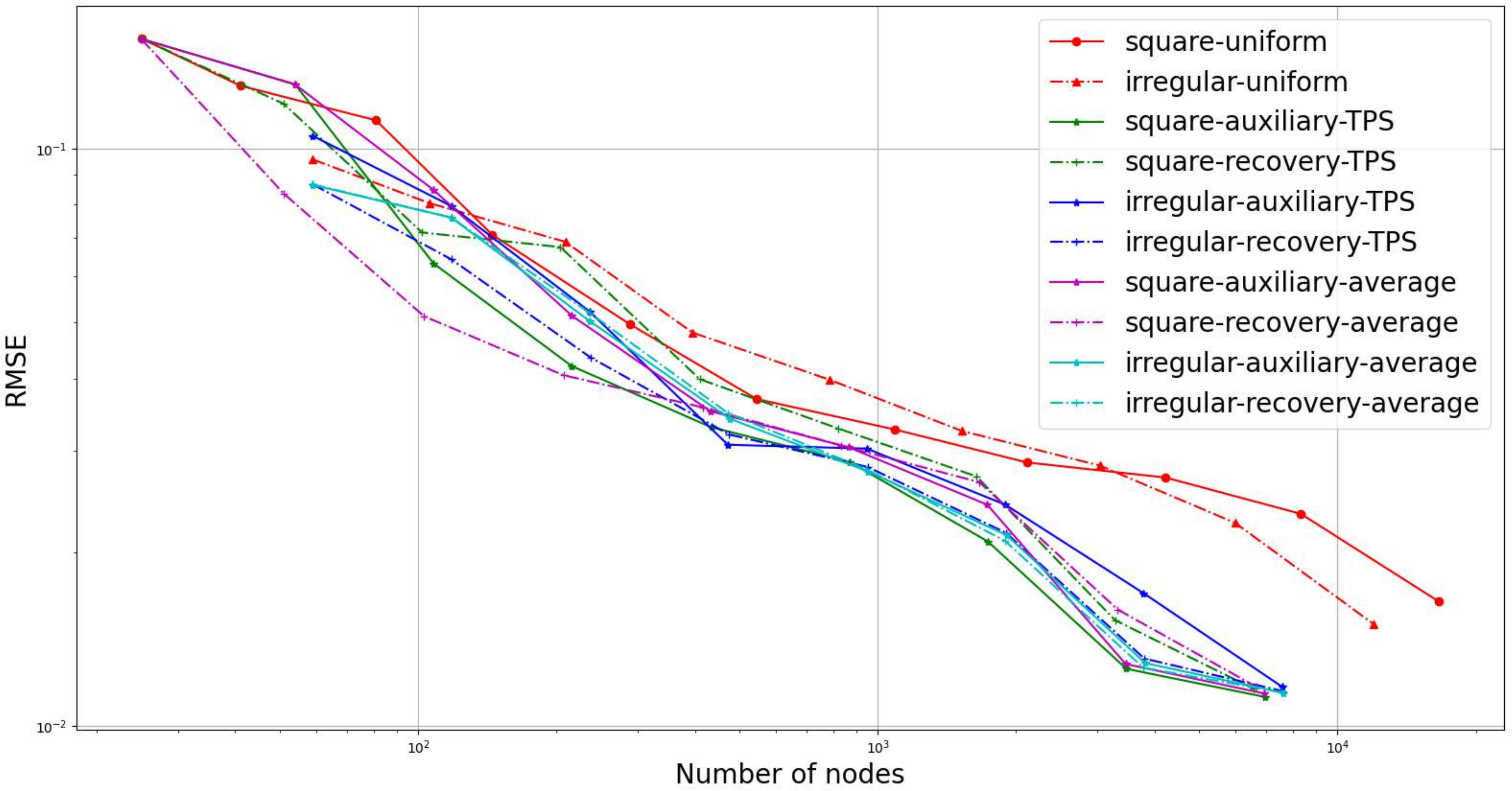}
		\caption{RMSE of TPSFEM for Mountain data.}
		\label{fig:convergence_7}
	\end{figure}

The convergence of the RMSE of the TPSFEM for the Mountain data is shown in Fig.~\ref{fig:convergence_7} and metrics of the final meshes are listed in Table~\ref{tab:metric_7}. The convergence curves correspond to TPSFEM in square or irregular domains with meshes refined uniformly or adaptively using the auxiliary problem or recovery-based error indicators. If values of new boundary nodes are approximated using the TPS during adaptive refinement, the curves are marked as TPS. Otherwise, if they are set as the average of their neighbouring nodes, the curves are marked as average. 

Adaptive refinement outperforms uniform refinement for the Mountain data in both square and irregular domains. Since most interpolated surfaces are smooth, the TPSFEM can achieve high accuracy with relatively coarse elements. The two error indicators have similar performance. Without adaptive refinement, the TPSFEM in irregular domains achieves lower RMSE with fewer nodes than square domains. The shorter distance between data and boundaries reduces the number of elements that do not cover any data point, which leads to higher efficiency. 

	\begin{table}
		\centering
		\caption{Metrics of TPSFEM for Mountain data}
		\label{tab:metric_7}
		\begin{tabular}{llllllll}
		\hline\noalign{\smallskip}
		Domain & Metric &  & \multicolumn{2}{c}{TPS approximation} & \multicolumn{2}{c}{Average} \\
		 &  &  Uniform & Auxiliary & Recovery & Auxiliary & Recovery \\
		\noalign{\smallskip}\hline\noalign{\smallskip}
		Square & \# nodes &  16641 & 6948 & 6564 & 6923 & 6656 \\
		 & Time &  0.504 & 0.225 & 0.226 & 0.299 & 0.216 \\
		 & RMSE &   0.0164 & 0.0112 & 0.0116 & 0.0114 & 0.0117 \\
		 & MAX &   0.45 & 0.353 & 0.356 & 0.356 & 0.355 \\
		\noalign{\smallskip}\hline\noalign{\smallskip}
		Irregular & \# nodes  & 11969 & 7576 & 7600 & 7634 & 7560 \\
		 & Time & 0.353 & 0.247 & 0.245 & 0.268 & 0.232 \\
		 & RMSE & 0.0150 & 0.0117 & 0.0115 & 0.0114 & 0.0114  \\
		 & MAX  & 0.441 & 0.356 & 0.352 & 0.355 & 0.354 \\
		\noalign{\smallskip}\hline
		\end{tabular}
	\end{table}



\begin{figure} [h]
		\centering
		\begin{subfigure}[b]{0.49\textwidth}
                \centering
			\includegraphics[scale=0.4]{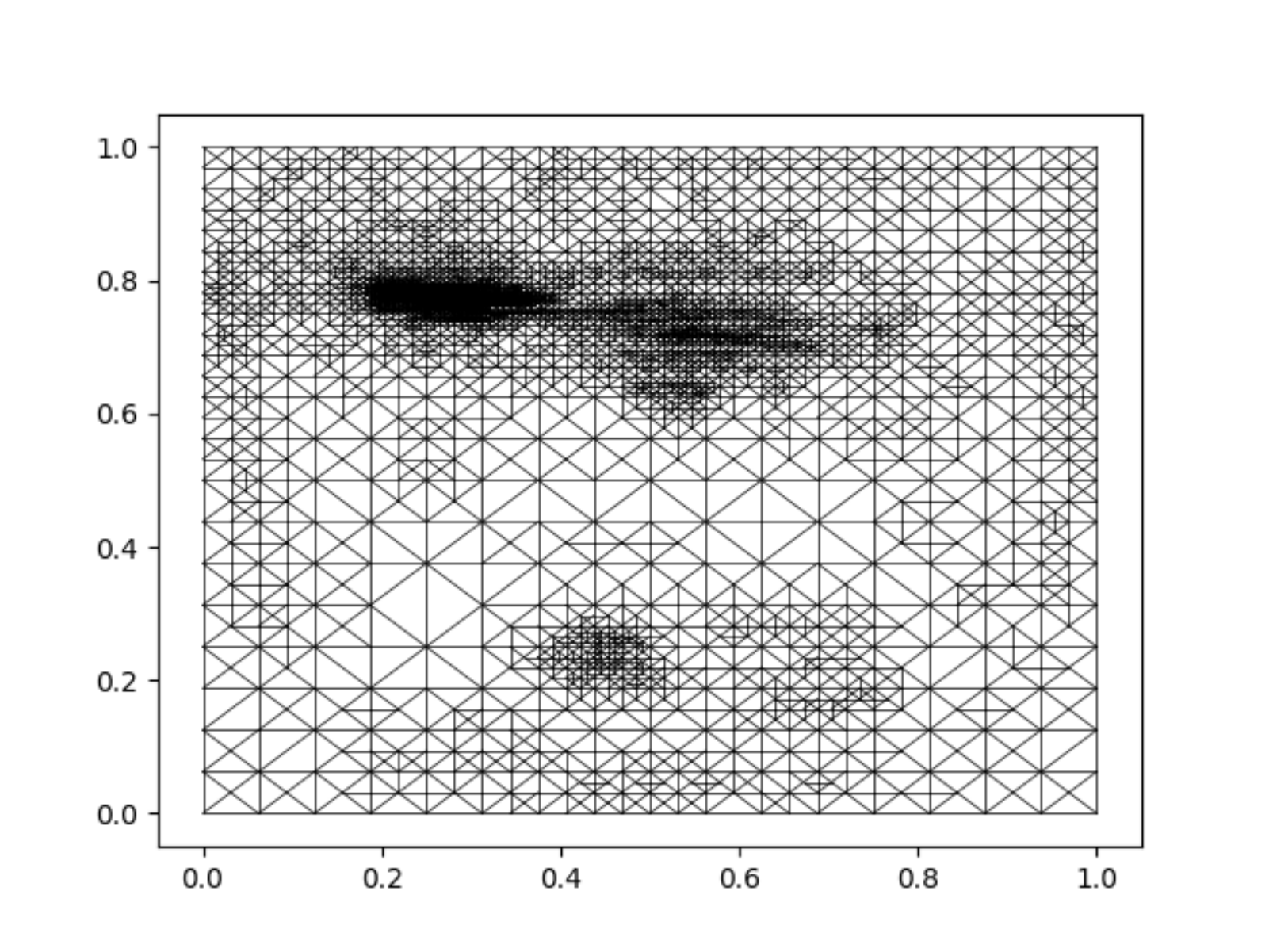}
			\caption{}
		\end{subfigure}
		\begin{subfigure}[b]{0.49\textwidth}
                \centering
			\includegraphics[scale=0.4]{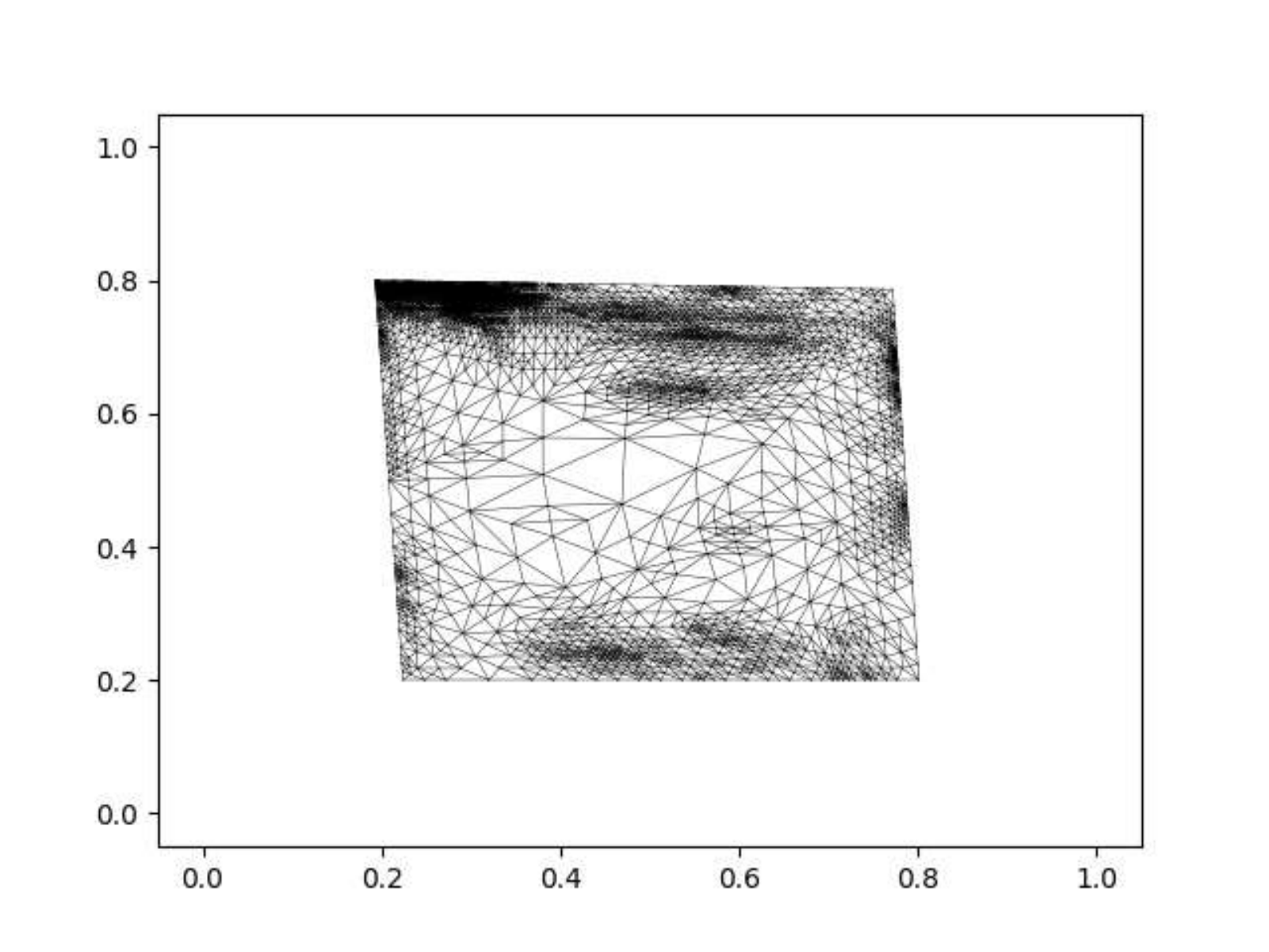}
			\caption{}
		\end{subfigure}
		\caption{FE meshes in (a) square domain; and (b) irregular domain. New boundary node values are approximated using TPS.}
		\label{fig:mesh_7_square_irregular_func}
	\end{figure}

When values of boundary nodes are approximated using the TPS, adaptive refinement achieves lower efficiency in irregular domains than in square domains. Adaptive meshes in irregular domains using the recovery-based error indicator achieve similar RMSE with $15.78\%$ more nodes than in square domains. The mesh refined using the auxiliary problem has higher RMSE and more nodes, which is markedly more inefficient. Two FE meshes generated using the recovery-based error indicator are shown in Fig.~\ref{fig:mesh_7_square_irregular_func}. The FE mesh in a square domain has $0\%$ of interior nodes near boundaries. In contrast, the adaptive mesh in the irregular domain in Fig.~\ref{fig:mesh_7_square_irregular_func}(b) has fine elements near boundaries and about $8.32\%$ of interior nodes close to boundaries. 

	\begin{figure} [h]
		\centering
		\begin{subfigure}[b]{0.49\textwidth}
                \centering
			\includegraphics[scale=0.4]{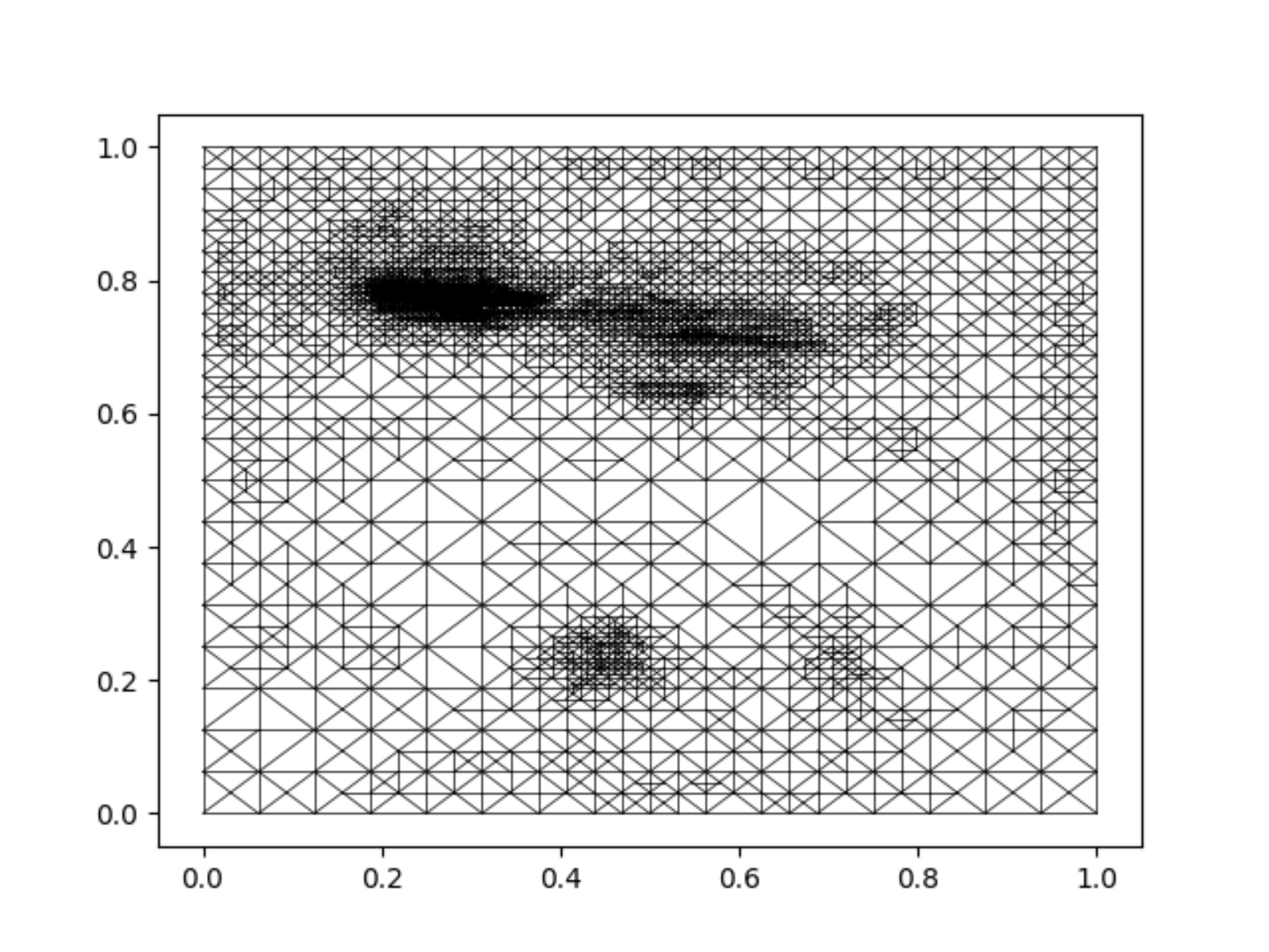}
			\caption{}
		\end{subfigure}
		\begin{subfigure}[b]{0.49\textwidth}
                \centering
			\includegraphics[scale=0.4]{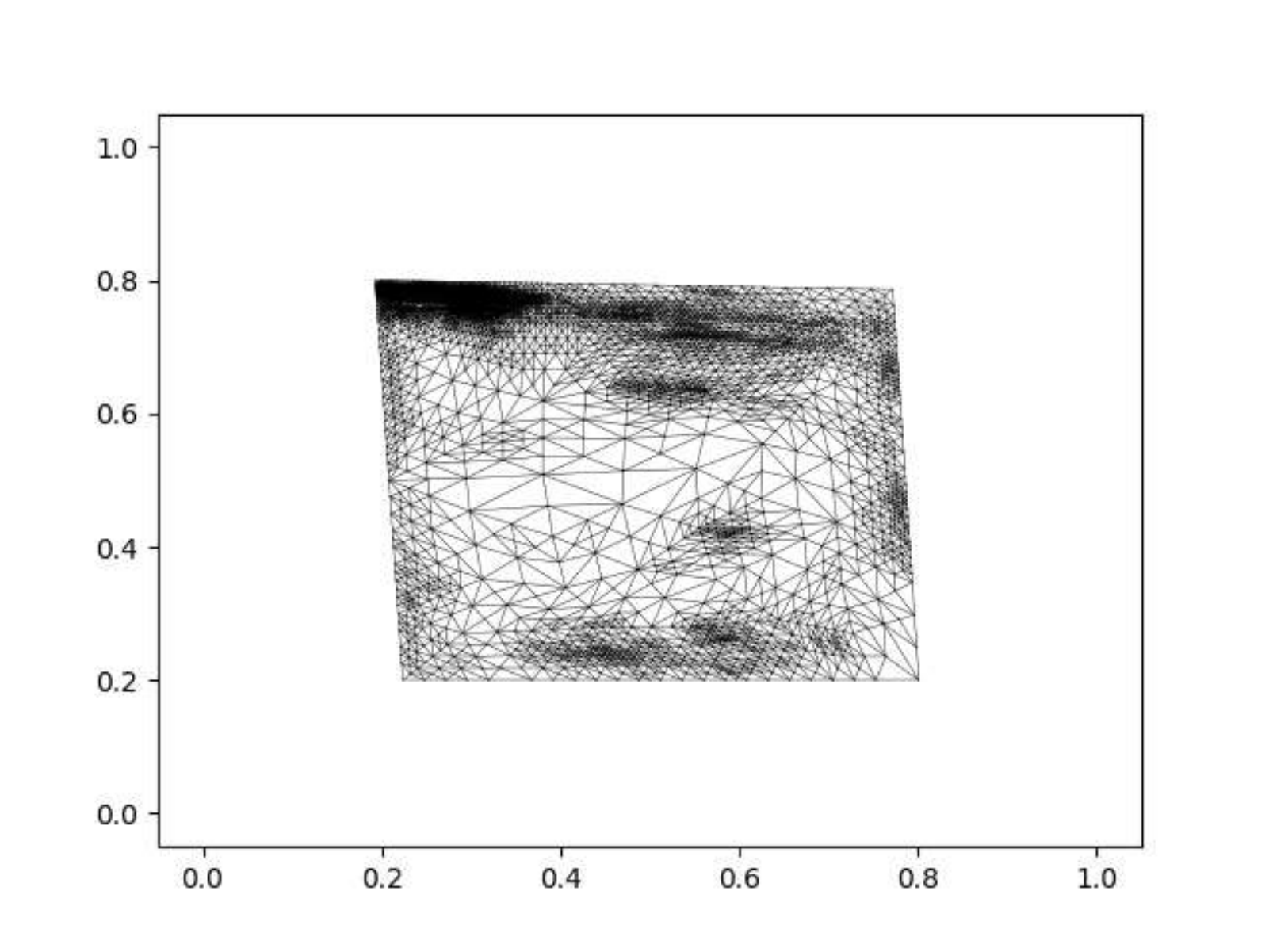}
			\caption{}
		\end{subfigure}
		\caption{FE meshes in (a) square domain; and (b) irregular domain. New boundary node values are set a average of neighbouring nodes.}
		\label{fig:mesh_7_square_irregular_average}
	\end{figure}

When new boundary node values during adaptive refinement are set using nodal average, over-refinement is markedly reduced as shown in Fig.~\ref{fig:mesh_7_square_irregular_average}. While no significant difference is observed between FE meshes in square domains as shown in Figs.~\ref{fig:mesh_7_square_irregular_func}(a) and~\ref{fig:mesh_7_square_irregular_average}(a), fewer fine elements near boundaries are found in Fig.~\ref{fig:mesh_7_square_irregular_average}(b) compared to Fig.~\ref{fig:mesh_7_square_irregular_func}(b). The ratio of interior nodes near boundaries is reduced to about 2.22\% in Fig.~\ref{fig:mesh_7_square_irregular_average}(b) and the accuracy of the solution like the peak at $[0.6,0.4]$ is improved with finer elements. Note that this peak is only identified by adaptive refinement in irregular domains and is not found in Figs.~\ref{fig:mesh_7_square_irregular_func}(a) and~\ref{fig:mesh_7_square_irregular_average}(a). 
Another interesting finding is that the auxiliary problem error indicator achieves the lowest RMSE in square domains with new boundary node values approximated using the TPS. Adaptive refinement identifies regions that do not require fine FEs and obtains similar efficiency.



\subsection{Canyon data}\label{sec:canyon}

\begin{figure} [h]
		\centering
		\includegraphics[scale=0.35]{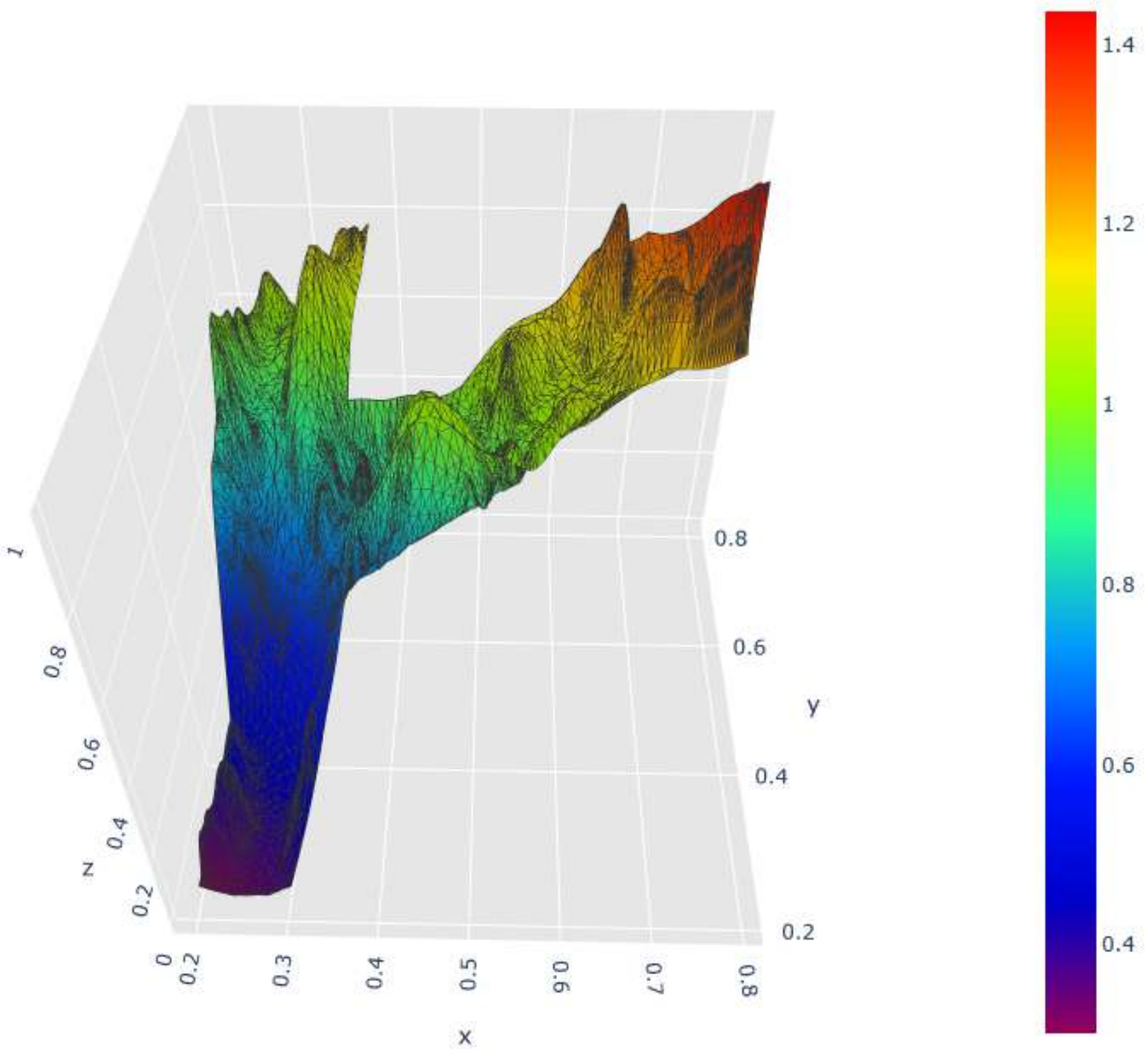}
		\caption{Triangulated surface of TPSFEM for Canyon data in an irregular domain.}
		\label{fig:trisurf_canyon}
	\end{figure}

A triangulated surface plot for the Canyon data generated using the auxiliary problem error indicator in an irregular domain is shown in Fig.~\ref{fig:trisurf_canyon}. The Canyon data consists of steep sides and a gradual bottom of the observed canyon. The steep sides are near the exterior of the Canyon data and these regions require finer elements to model. In contrast, the bottom of the canyon and regions without data points only need coarse elements.
	
\begin{figure} [h]
		\centering
		\includegraphics[scale=0.26]{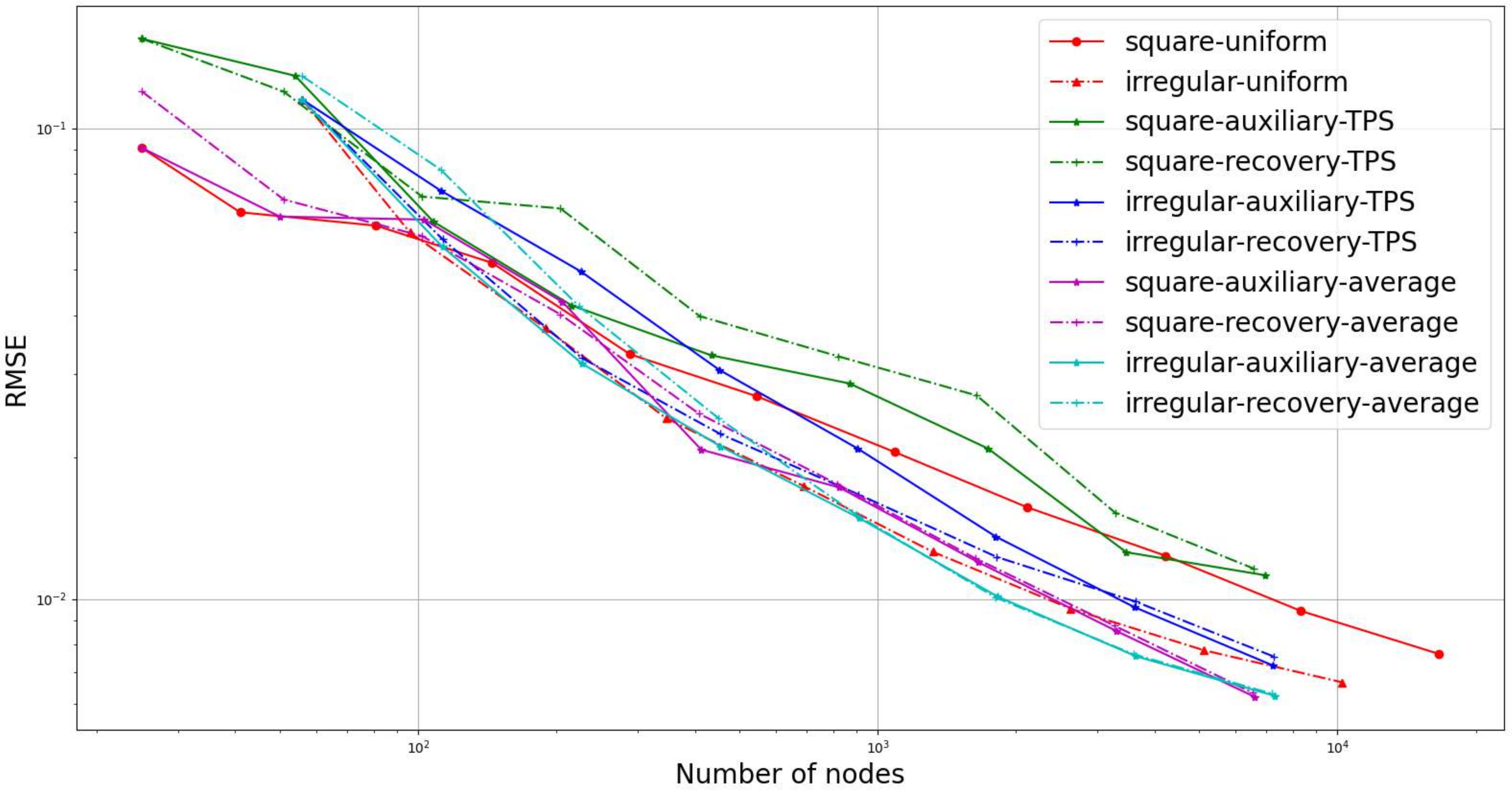}
		\caption{RMSE of TPSFEM for Canyon data.}
		\label{fig:convergence_9}
	\end{figure}

The convergence of the RMSE of the TPSFEM for the Canyon data is shown in Fig.~\ref{fig:convergence_9} and the metrics of the final meshes are shown in Table~\ref{tab:metric_9}. It is shown that adaptive meshes significantly improve the efficiency of the Canyon data in square domains with new boundary node values calculated by the TPS. For example, the adaptive mesh of the auxiliary problem error indicator achieves about $15.74\%$ of RMSE using $39.74\%$ number of nodes compared to the uniform mesh. Adaptive refinement produces finer elements on steep sides and coarser elements on regions without data. This allows the TPSFEM to achieve high RMSE while retaining efficiency. The auxiliary problem error indicator outperforms the recovery-based error indicator in all cases of Table~\ref{tab:metric_9} as it better detects oscillatory behaviour using local data points.

	\begin{table}
		\centering
		\caption{Metrics of TPSFEM for Canyon data}
		\label{tab:metric_9}
		\begin{tabular}{lllllll}
		\hline\noalign{\smallskip}
		Domain & Metric &  & \multicolumn{2}{c}{TPS approximation} & \multicolumn{2}{c}{Average} \\
		 &  & Uniform & Auxiliary & Recovery & Auxiliary & Recovery \\
		\noalign{\smallskip}\hline\noalign{\smallskip}
		Square & \# nodes  & 16641 & 6612 & 6528 & 6608 & 6534 \\
		 & Time & 0.448 & 0.231 & 0.223 & 0.228 & 0.235 \\
		 & RMSE & 0.0394 & 0.00620 & 0.00665 & 0.00619 & 0.00632 \\
		 & MAX & 0.0524 & 0.0623 & 0.0650 & 0.0656 & 0.0605 \\
		\noalign{\smallskip}\hline\noalign{\smallskip}
		Irregular & \# nodes & 10241 & 7238 & 7265 & 7294 & 7204 \\
		 & Time & 0.299 & 0.236 & 0.242 & 0.232 & 0.229 \\
		 & RMSE  & 0.00665 & 0.00722 & 0.00754 & 0.00623 & 0.00629  \\
		 & MAX  & 0.118 & 0.117 & 0.116 & 0.116 & 0.117 \\
		\noalign{\smallskip}\hline
		\end{tabular}
	\end{table}

\begin{figure} [h]
		\centering
		\begin{subfigure}[b]{0.49\textwidth}
                \centering
			\includegraphics[scale=0.4]{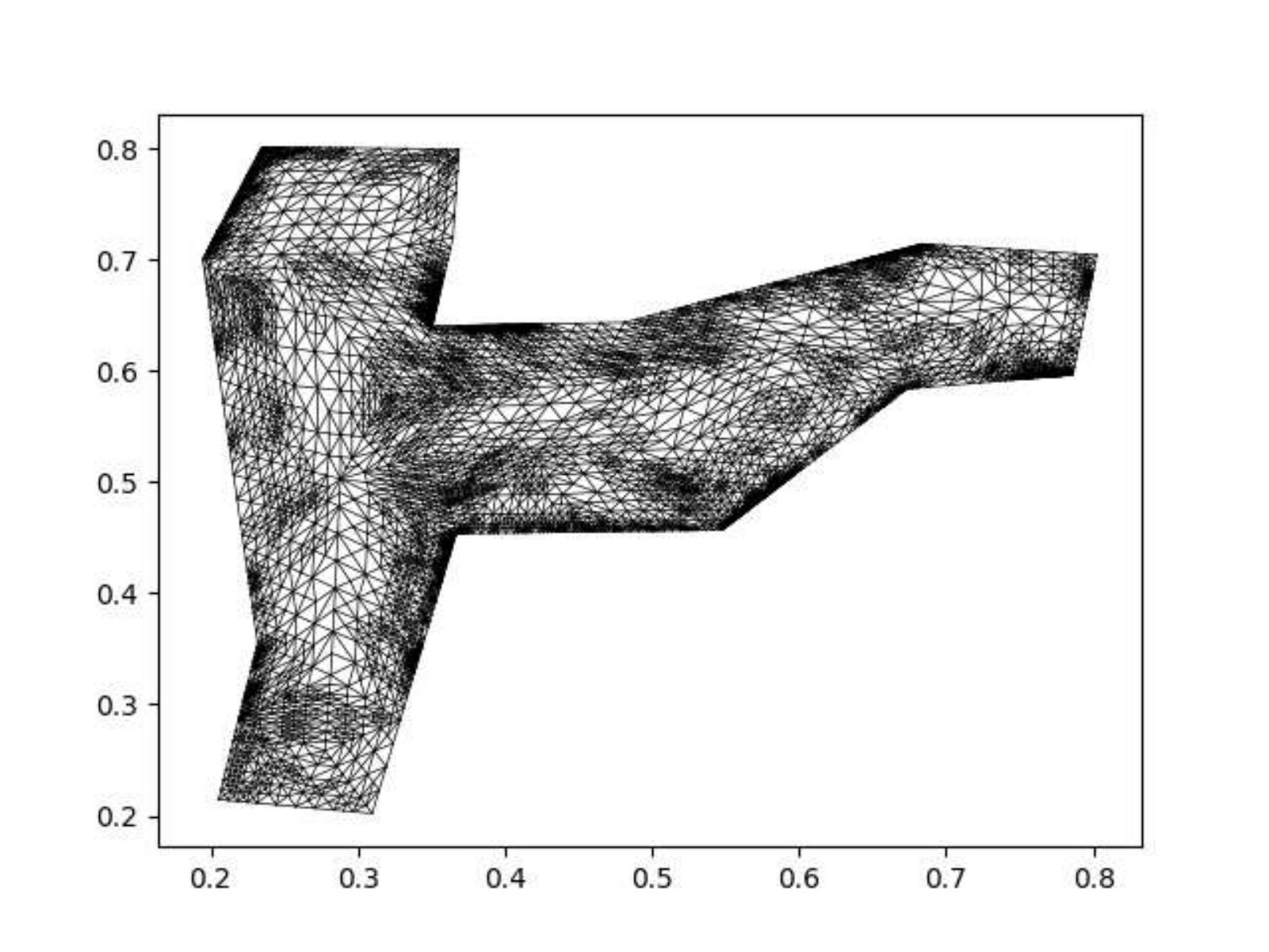}
			\caption{}
		\end{subfigure}
		\begin{subfigure}[b]{0.49\textwidth}
                \centering
			\includegraphics[scale=0.4]{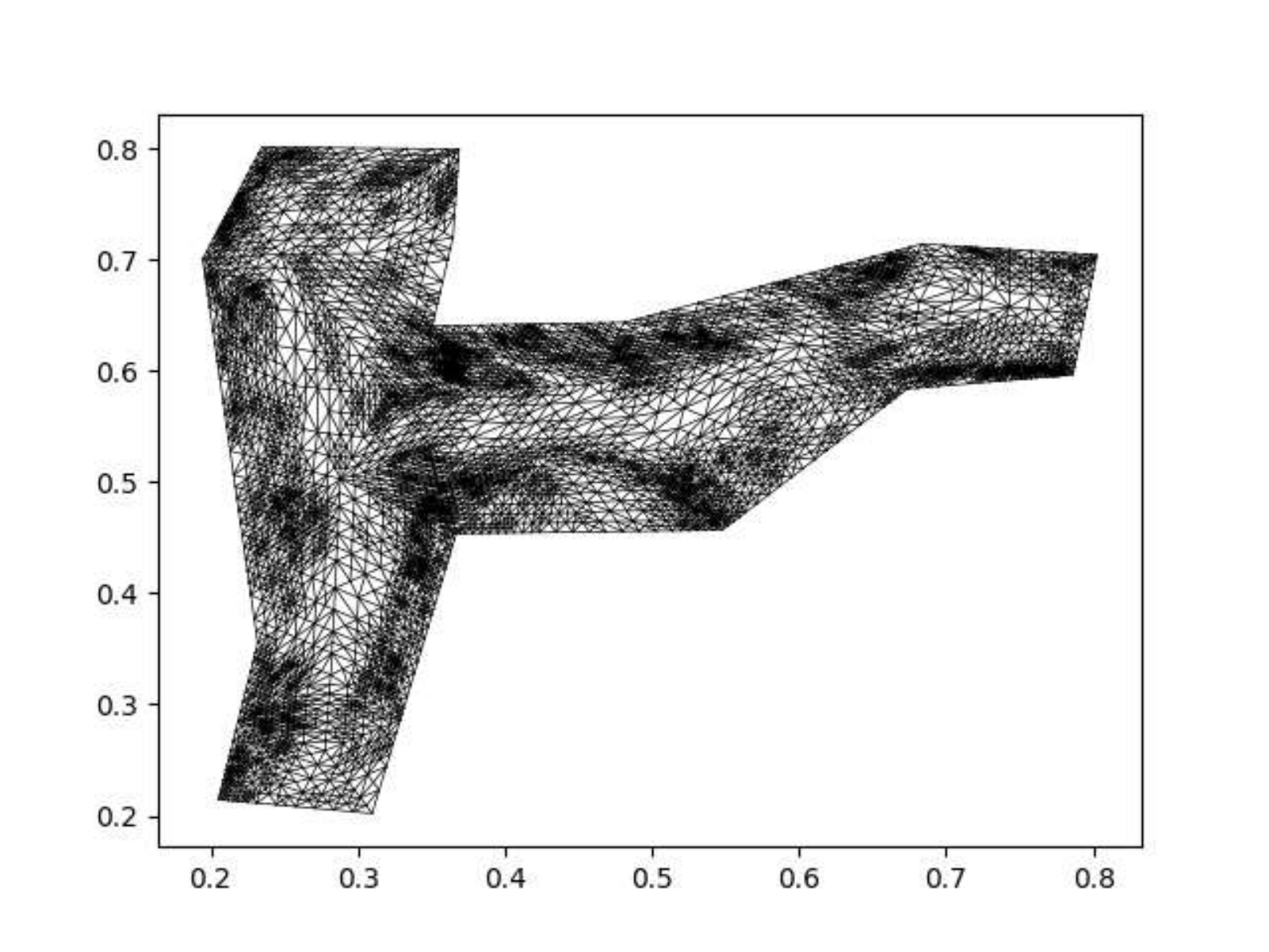}
			\caption{}
		\end{subfigure}
		\caption{FE meshes in irregular domain generated with new boundary node values (a) approximated using TPS; (b) set as average of neighbouring nodes.}
		\label{fig:mesh_9_irregular_func_average}
	\end{figure}

Since the Canyon data only covers a small portion of the domain as shown in Fig.~\ref{fig:data_dist}, irregular domains significantly improve the TPSFEM's efficiency. However, adaptive refinement fails to improve the efficiency in irregular domains when new boundary node values are initialised using the TPS approximations. An adaptive mesh generated using the auxiliary problem error indicator is shown in Fig.~\ref{fig:mesh_9_irregular_func_average}(a), where severe over-refinement is observed near boundaries and about 33.43\% of interior nodes are close to boundaries. Consequently, the RMSE for both the error indicators in irregular domains has significantly higher RMSE and numbers of nodes than those in square domains. Over-refinement near boundaries has a major impact on the Canyon data compared to the Mountain data.

After adopting the nodal average approach, over-refinement near boundaries is markedly reduced as shown in Fig.~\ref{fig:mesh_9_irregular_func_average}(b) and the ratio of interior nodes near boundaries is reduced to about 5.61\%. Consequently, the RMSE of final adaptive meshes in both square and irregular domains achieves similar RMSE as those in square domains. In addition, the TPSFEM has much higher MAX in irregular domains than in square domains. Since the Canyon data has more oscillatory surfaces near boundaries than the other two data sets, the approximated boundary conditions have larger impacts on $s$, especially near boundaries. 



\subsection{River data}\label{sec:river}

\begin{figure} [h]
		\centering
		\includegraphics[scale=0.35]{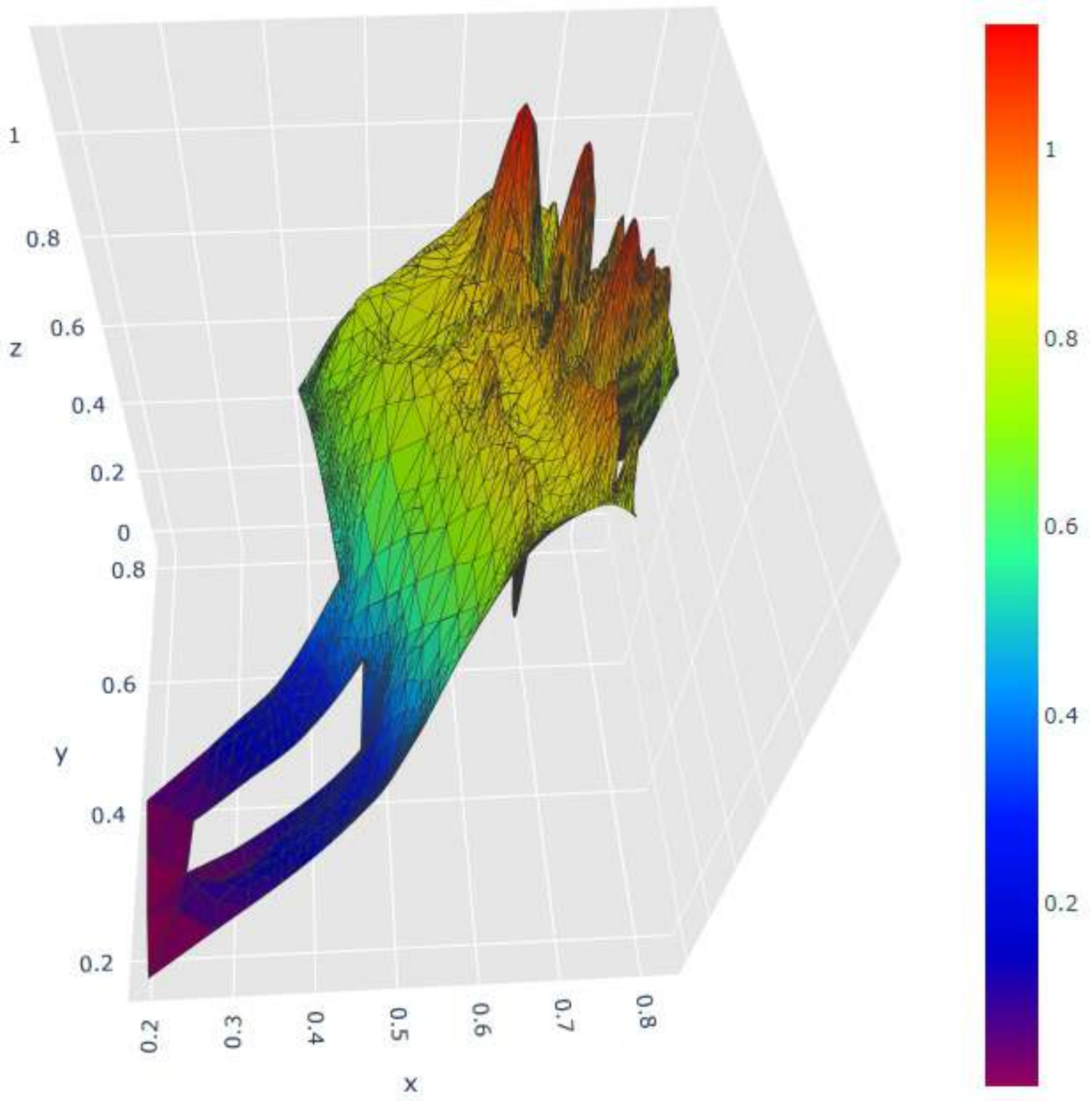}
		\caption{Triangulated surface of TPSFEM for River data in an irregular domain.}
		\label{fig:trisurf_river}
	\end{figure}

A triangulated surface plot for the River data generated using the recovery-based error indicator in its irregular domain is shown in Fig.~\ref{fig:trisurf_river}. The River data is the most complicated data set with various routes and data densities as described in Sect.~\ref{sec:data}. The sampled data concentrates at the river mouth, where currents change and river beds are eroded by the flow. This leads to oscillatory surfaces and deep holes under the river mouth near point $[0.8,0.6]$. The other regions are relatively smooth and only require coarse elements. Only elements near the river mouth are refined with finer elements.

\begin{figure} [h]
		\centering
		\includegraphics[scale=0.26]{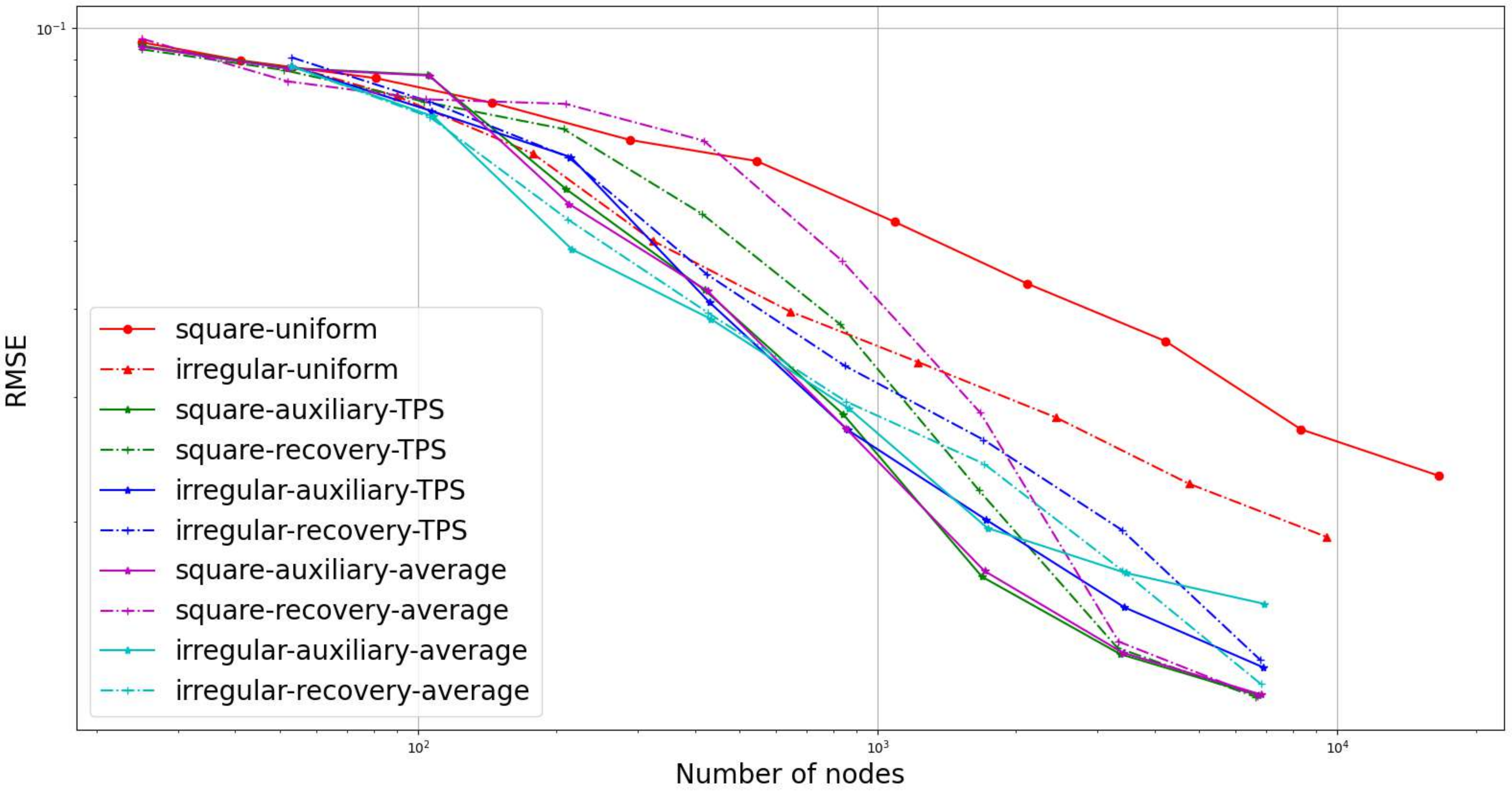}
		\caption{RMSE of TPSFEM for River data.}
		\label{fig:convergence_10}
	\end{figure}

The convergence of the RMSE of the TPSFEM for the River data is shown in Fig.~\ref{fig:convergence_10} and metrics of the final meshes are shown in Table~\ref{tab:metric_10}. Adaptive meshes in square domains achieve much higher convergence rates compared to uniform meshes for the River data. Since the river mouth region is significantly more oscillatory than other regions, adaptive meshes are more efficient than uniform meshes similar to the Canyon data. Both the error indicators have similar performance in this data set.

	\begin{table}
		\centering
		\caption{Metrics of TPSFEM for River data}
		\label{tab:metric_10}
		\begin{tabular}{lllllll}
		\hline\noalign{\smallskip}
		Domain & Metric & & \multicolumn{2}{c}{TPS approximation} & \multicolumn{2}{c}{Average} \\
		 &  & Uniform & Auxiliary & Recovery & Auxiliary & Recovery \\
		\noalign{\smallskip}\hline\noalign{\smallskip}
		Square & \# nodes  & 16641 & 6740 & 6642 & 6828 & 6688 \\
		 & Time  & 0.465 & 0.234 & 0.228 & 0.227 & 0.248 \\
		 & RMSE  & 0.0232 & 0.0113 & 0.0113 & 0.0114 & 0.0113  \\
		 & MAX  & 0.379 & 0.194 & 0.193 & 0.199 & 0.193 \\
		\noalign{\smallskip}\hline\noalign{\smallskip}
		Irregular & \# nodes  & 9488  & 6928 & 6792 & 6880 & 6816 \\
		 & Time & 0.304 & 0.238 & 0.243 & 0.230 & 0.240 \\
		 & RMSE  & 0.0190 & 0.0152 & 0.0127 & 0.0124 & 0.0118 \\
		 & MAX  & 0.283 & 0.204 & 0.174 & 0.172 & 0.172 \\
		\noalign{\smallskip}\hline
		\end{tabular}
	\end{table}
	
\begin{figure} [h]
		\centering
		\begin{subfigure}[b]{0.49\textwidth}
                \centering
			\includegraphics[scale=0.4]{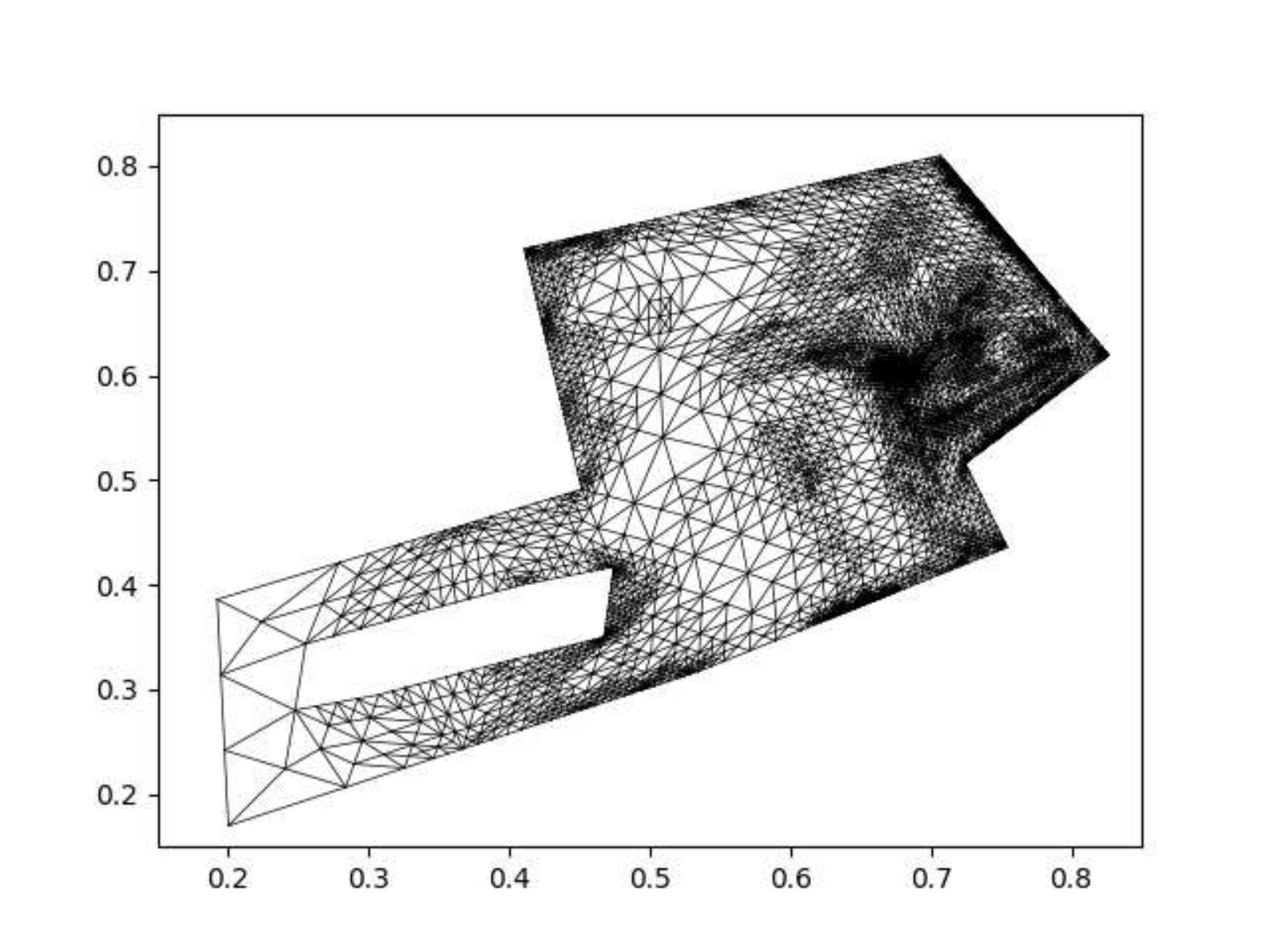}
			\caption{}
		\end{subfigure}
		\begin{subfigure}[b]{0.49\textwidth}
                \centering
			\includegraphics[scale=0.4]{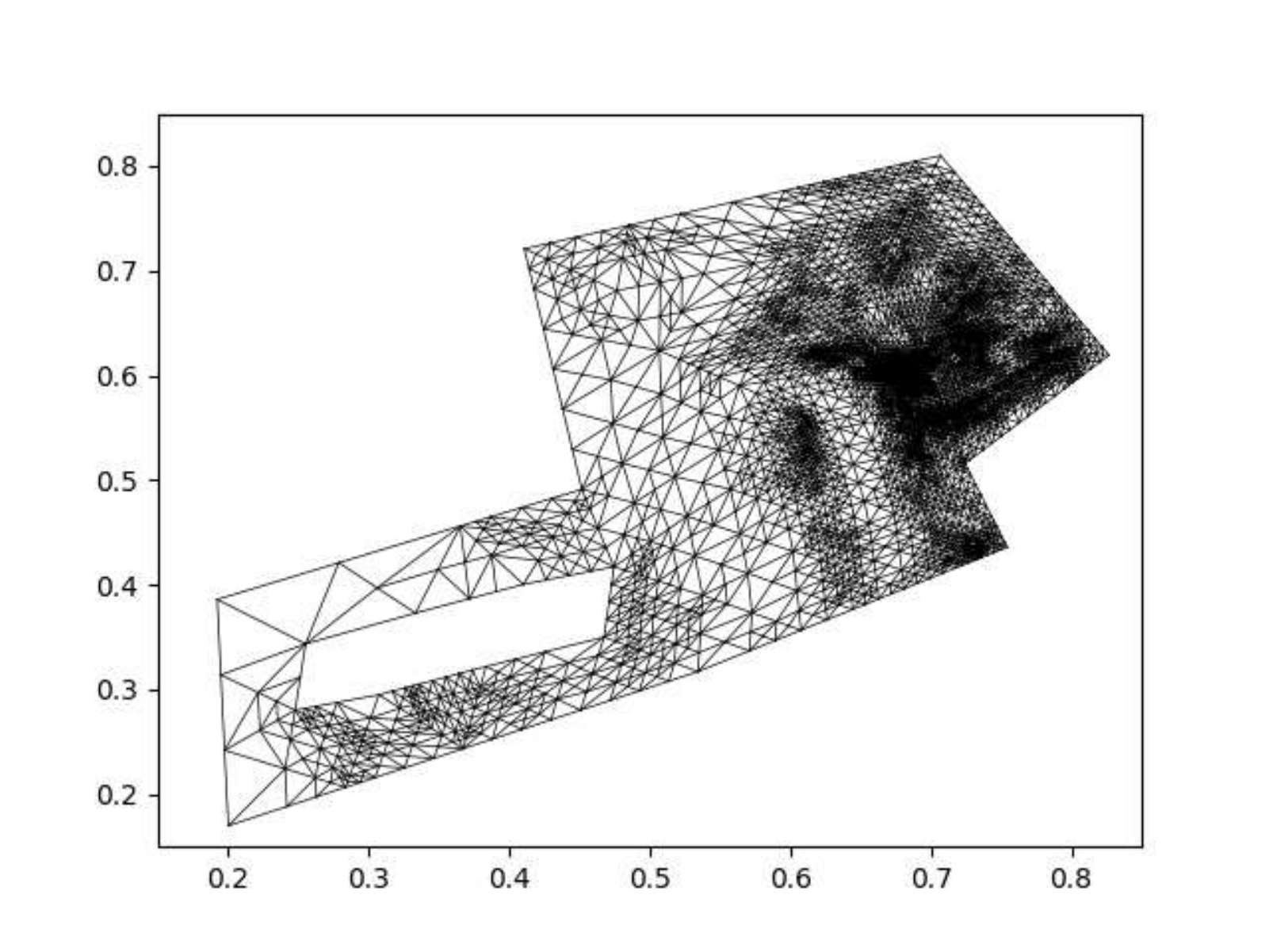}
			\caption{}
		\end{subfigure}
		\caption{FE meshes in irregular domain generated with values of new boundary nodes (a) approximated using TPS; (b) set as average of neighbouring nodes.}
		\label{fig:mesh_10_irregular_func_average}
	\end{figure}

The TPSFEM in irregular domains with uniform refinement obtains $81.90\%$ of RMSE using $57.02\%$ of the number of nodes compared to that in square domains. The improvement is not as significant as adaptive refinement, which achieves $48.88\%$ of RMSE using $40.50\%$ of the number of nodes in square domains. Adaptive refinement also improves the MAX in all test cases. In addition, adaptive refinement using the TPS approximations leads to higher RMSE and number of nodes in irregular domains compared to square domains. This is caused by small oscillatory regions surrounding the river mouth. Since the area is close to the boundaries of the irregular domain, the error indicators are misled and over-refined in multiple regions near the boundaries of the FE mesh as illustrated in Fig.~\ref{fig:mesh_10_irregular_func_average}(a). Over-refinement is significantly reduced using the nodal average approach as shown in Fig.~\ref{fig:mesh_10_irregular_func_average}(b), where only a small amount of over-refinement is observed. The ratio of interior nodes near boundaries has also been drastically reduced from about $13.18\%$ to $1.19\%$. 

	
\subsection{Comparison to RBFs}\label{sec:csrbf}

The performance of the TPSFEM is compared to the TPS and two CSRBFs from~\cite{wendland1995piecewise} and~\cite{buhmann1998radial} to illustrate its effectiveness. Let $\bar{\bm{x}}_i$ be a control point, where $1 \le i \le \bar{n}$. Its kernel is defined as $\phi_{i}(r)$ in Table~\ref{tab:comp_rbf}, where~$r=||\bar{\bm{x}}_i-\bm{x}||_{2}$ is the Euclidean distance between $\bar{\bm{x}}_i$ and $\bm{x}$. The TPS has global support and the two CSRBFs have compact support as their kernels $\phi(r)$ are defined as $\phi(r)=0$ if $r>1$. The kernels of the two CSRBFs are scaled using a radius of support~$\rho$ and their values are computed as $\phi(r/\rho)$, which markedly reduces the computational costs. In the experiments, two values of~$\rho$ are chosen for each CSRBF such that the compact support of each $\bar{\bm{x}}_i$ covers a fixed number of data points~\citep{deparis2014rescaled}. For example, the CSRBFs for the Mountain data with a radius of support 0.0691 and 0.120 cover approximately 100 or 300 sampled points, respectively.

	\begin{table}
		\centering
		\caption{Kernels of TPS and CSRBFs.}
		\label{tab:comp_rbf}
		\begin{tabular}{ll}
			\hline\noalign{\smallskip}
			Technique & Kernel $\phi$ \\
			\noalign{\smallskip}\hline\noalign{\smallskip}
			TPS & $r^{2}\log(r)$ \\	
			Buhmann & $1/15+19/6r^{2}-16/3r^{3}+3r^{4}-16/15r^{5}+1/6r^{6}+2r^{2}\log(r)$ \\			
			Wendland & $(1-r)^{4}(4r+1)$  \\
			\noalign{\smallskip}\hline
		\end{tabular}
	\end{table}

The control points of the RBFs are sampled by finding data points close to nodes of a uniform rectangular mesh over the domain. Given a mesh with mesh size $\bar{h}$, regions where the distance between the nearest nodes and data point is greater than $\bar h/3$ are ignored. Thus, the maximum and minimum distance between neighbouring control points is bounded. The smoothing parameters of the RBFs are calculated using the GCV.

The efficiency of the RBFs and TPSFEM are compared using statistics listed in Table~\ref{tab:metric_rbf}. 
The RBFs have $\bar{n}$ basis and the TPSFEM has $m$ nodes and one basis function for each node. The computational costs are evaluated using the ratio and numbers of nonzero entries in the system of equations and the time to solve the system, which is measured in seconds. The systems of the TPS are solved using a direct solver and the systems of the CSRBFs and TPSFEM are solved using the sparse solver (PyPardiso). The TPSFEM generated in irregular domains using the recovery-based error indicator and nodal average for the three data sets have also been listed in Table~\ref{tab:metric_rbf} for comparison.

	\begin{table}
		\centering
		\caption{Metrics of CSRBF for Mountain, Canyon and River data.}
		\label{tab:metric_rbf}
        \resizebox{\columnwidth}{!}{
		\begin{tabular}{llllllll}
		\hline\noalign{\smallskip}
		Data & Metric & TPS & \multicolumn{2}{c}{Buhmann} & \multicolumn{2}{c}{Wendland} & TPSFEM \\
		\noalign{\smallskip}\hline\noalign{\smallskip}
		Mountain & Radius & & 0.0691 & 0.120 & 0.0691 & 0.120  &  \\
		 & \# basis & 6673 & 6673 & 6673 & 6673 & 6673 &  7560 \\
		 & \# nonzero & 44,528,929 & 1,809,811 & 4,959,647 & 1,809,811 & 4,959,64 &  426,448 \\
		 & Ratio & 100\% & 4.06\% & 11.1\% & 4.06\% & 11.1\% & 0.044\% \\
		 & Time  & 2.56 & 0.777 & 1.38 & 0.779 & 1.45 & 0.232 \\
		 & RMSE  & 0.0166 & 0.0165 & 0.0165 & 0.0171 & 0.0172 & 0.0114  \\
		 & MAX  & 0.575 & 0.573 & 0.573 & 0.591 & 0.591 & 0.354 \\
		\noalign{\smallskip}\hline\noalign{\smallskip}
		Canyon & Radius & & 0.0684 & 0.118 & 0.0684 & 0.118  &  \\
		 & \# basis  &  6807 & 6807 & 6807 & 6807 & 6807 &  7204 \\
		 & \# nonzero & 46,335,249 & 5,984,357 & 10,450,523 & 5,984,357 & 10,450,523 & 409,826 \\
		 & Ratio & 100\% & 12.9\% & 22.6\% & 12.9\% & 22.6\% & 0.046\% \\
		 & Time & 3.03 & 1.48 & 2.35 & 1.42 & 2.41 & 0.229 \\
		 & RMSE  & 0.00658 & 0.00668 & 0.00658 & 0.00680 & 0.00679 &  0.00629  \\
		 & MAX  & 0.0808 & 0.0795 & 0.0798 & 0.0914 & 0.0904 &  0.117 \\
		\noalign{\smallskip}\hline\noalign{\smallskip}
		River & Radius & & 0.0699 & 0.121 & 0.0699 & 0.121  &  \\
		 & \# basis &  6517 & 6517 & 6517 & 6517 & 6517 &  6816 \\
		 & \# nonzero & 42,471,289 & 4,240,983 & 9,815,249 & 4,240,983 & 9,815,249 & 399,667 \\
		 & Ratio & 100\% & 9.97\% & 23.2\% & 9.97\% & 23.2\% & 0.05\% \\
		 & Time & 2.07 & 1.16 & 2.36 & 1.15 & 2.35 & 0.240 \\
		 & RMSE  & 0.0144 & 0.0146 & 0.0145 & 0.0139 & 0.0140 &  0.0118  \\
		 & MAX  & 0.370 & 0.373 & 0.372 & 0.344 & 0.342 &  0.172 \\
		\noalign{\smallskip}\hline
		\end{tabular}
          }
	\end{table}

As shown in Table~\ref{tab:metric_rbf}, the CSRBFs from Buhmann and Wendland have similar RMSE and MAX to the TPS for the three data sets. Since the systems of the CSRBFs are much sparser than the dense systems of the TPS as shown by the nonzero ratios, the CSRBFs require significantly less time to solve the system. The CSRBF of Buhmann achieves higher efficiency for the Mountain and Canyon data, but it underperforms for the River data compared to the CSRBF of Wendland. The TPSFEM with uniform meshes in square domains are markedly less efficient compared to the CSRBFs as shown in Tables~\ref{tab:metric_7},~\ref{tab:metric_9} and~\ref{tab:metric_10}. However, the TPSFEM with adaptive refinement achieves much lower RMSE and MAX for all three data sets. Additionally, their systems are much sparser and take less time to solve.

The performance of the CSRBFs is affected by factors like the radius of support~$\rho$. Table~\ref{tab:metric_rbf} contains attributes of the CSRBFs with~$\rho$ that cover about 100 and 300 sample points corresponding to columns 4, 6 and rows 5, 7, respectively. While the increase in $\rho$ may reduce the RMSE and MAX of the CSRBFs in other studies, it showed little improvement for the three data sets. Note that the performance of the CSRBFs can be improved using a localised radius of support, preconditioners or parallel solvers~\citep{deparis2014rescaled}. Nevertheless, this experiment demonstrates the TPSFEM's effectiveness.

\section{Conclusions}\label{sec:conclusions}

We study the TPSFEM and adaptive mesh refinement in irregular domains with complex geometric shapes, which better represent the geometry for the domain of interest. Unknown Dirichlet boundary conditions for real-world data sets are approximated using the TPS with a small subset of sampled data points. They improve the accuracy of the TPSFEM for irregular domains with data points close to boundaries. Over-refinement of FEs near boundaries weakens the performance of adaptive refinement for the TPSFEM. An alternative is proposed to initialise new boundary node values as the average of neighbouring nodes. This improves the stability of adaptive refinement when exact boundary conditions are unavailable.

Numerical experiments are conducted to investigate the performance of the TPSFEM and adaptive refinement in square and irregular domains.  The experiments use three geological surveys with varied distribution patterns and irregular domains. Two approaches to initialising new boundary node values during adaptive refinement are compared. The results show that adaptive refinement markedly improves the efficiency of the TPSFEM for all three surveys in both square and irregular domains. While over-refinement of FEs occurs in irregular domains, it is significantly reduced using the proposed averaging initialisation for new boundary nodes. Irregular domains also increase the efficiency of the TPSFEM with uniform refinement, though adaptive refinement achieves slightly higher convergence in square domains. The recovery-based error indicator is more stable and slightly overperforms the auxiliary problem error indicator. The comparison of the TPSFEM, TPS and two chosen CSRBFs shows that the TPSFEM achieves lower RMSE with costs. Overall, concerning efficiency and flexibility for digital terrain modelling in irregular domains, our findings show the TPSFEM is worthy of recommendation. Our future work will focus on extending it to three-dimensional data sets and improving its accuracy for surfaces with abrupt changes.

\begin{acknowledgements}
I would like to express my great appreciation to Dr. Linda Stals for her valuable and constructive suggestions during the planning and development of this research work. Her willingness to give her time so generously has been very much appreciated. This work was supported by the Natural Science Foundation of Xiamen [Grant No. 20237101] and the Scientific Research Funds of Huaqiao University.


    
\end{acknowledgements}

\section*{Data availability}
The data that support the findings of this study are openly available in Bitbucket at https://bitbucket.org/fanglishan/tpsfem-dtm-data/src/main/.


\bibliographystyle{MG}
{\footnotesize
\bibliography{tpsfem_dtm_fang}} 

\appendix

\section{Triangulation in irregular domain}\label{app:domain}

The triangulation in irregular domains is more complicated than in regular rectangular domains and many software packages are available. \cite{ramsay2002spline} used a PDE toolbox of MATLAB to build arbitrarily shaped triangulated domains. While every data point is a vertex, not every vertex corresponds to a data point. However, this approach only applies to problems with few data points.
	
\begin{figure} [h]
		\centering
		\begin{subfigure}[b]{0.49\textwidth}
                \centering
			\includegraphics[scale=0.4]{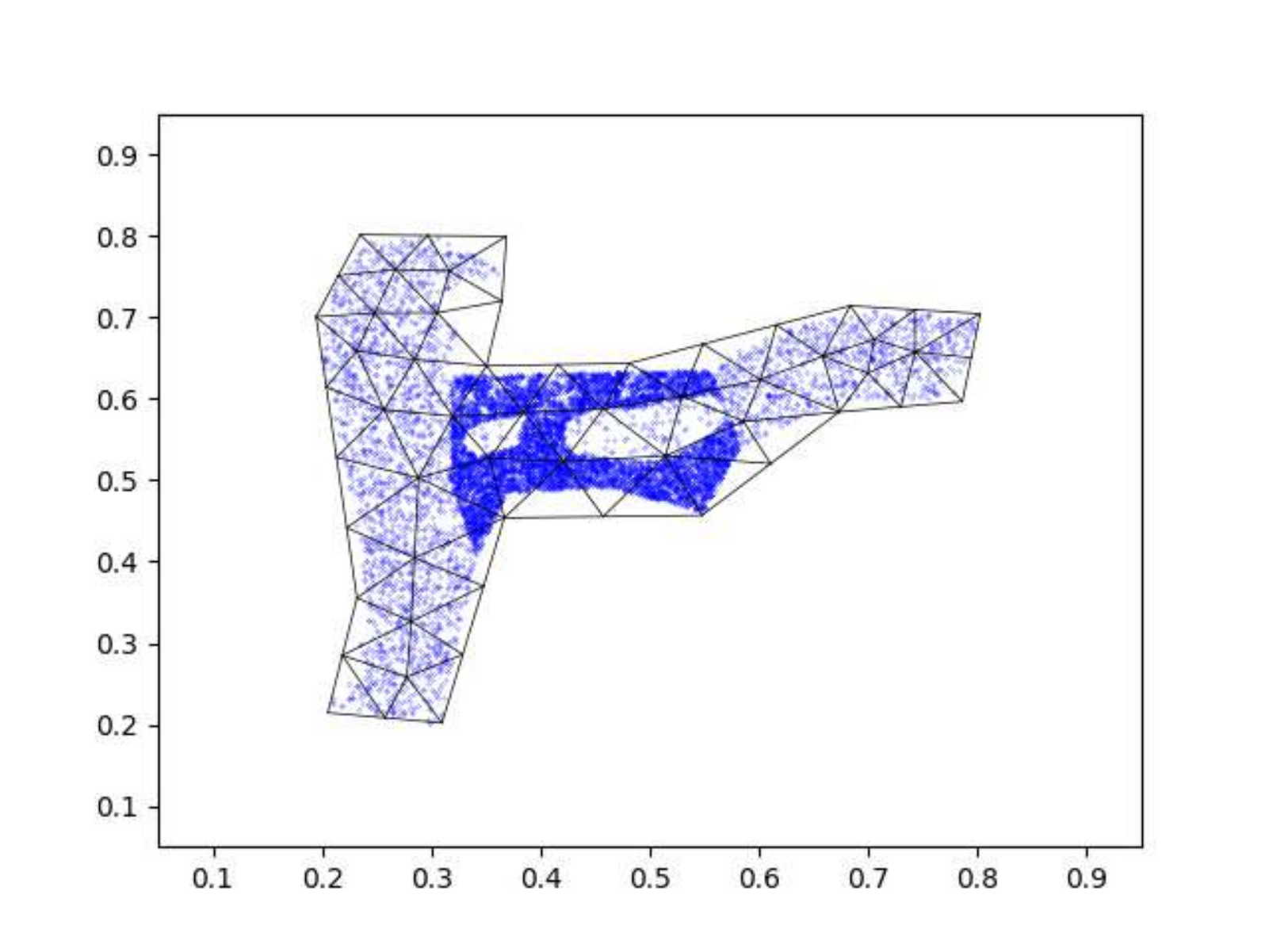}
			\caption{}
		\end{subfigure}
		\begin{subfigure}[b]{0.49\textwidth}
                \centering
			\includegraphics[scale=0.4]{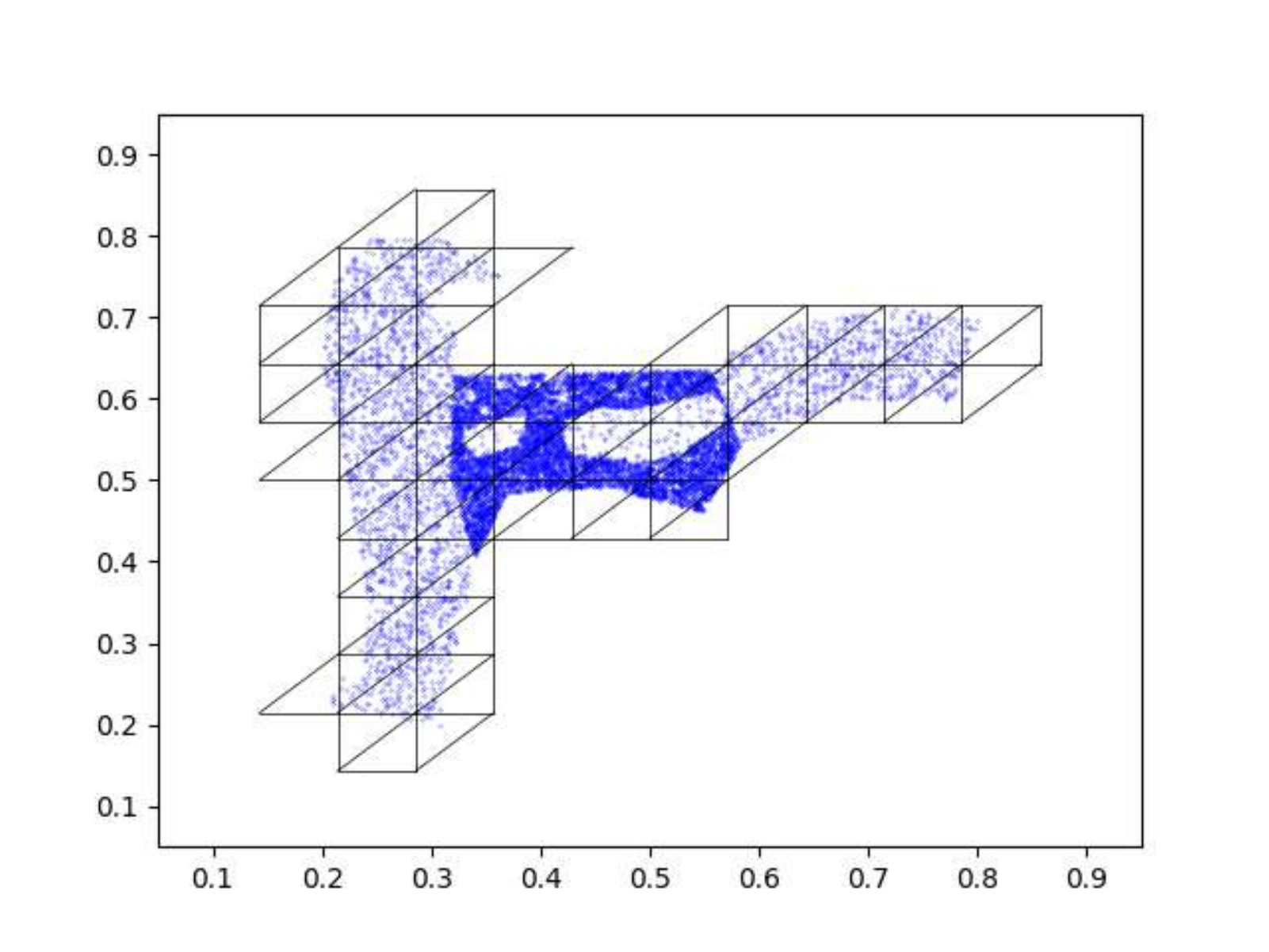}
			\caption{}
		\end{subfigure}
		\caption{FE meshes generated by (a) triangulation in an irregular domain; and (b) removing elements from square domain. The data points are represented as blue dots.}
		\label{fig:mesh_polygon_trim}
	\end{figure}

Irregular domains of the TPSFEM are built using Python packages alphashape\footnote{alphashape, https://pypi.org/project/alphashape/} and pygmsh\footnote{pygmsh, https://pypi.org/project/pygmsh/}. Firstly, a set of vertices are obtained to fit the exterior of the observed data using alphashape, which generalises bounding polygons containing the data~\citep{xu2003automatic}.  A triangulated mesh is then generated using pygmsh. An example FE mesh is shown in Fig.~\ref{fig:mesh_polygon_trim}(a), which fits the observed data closely. Another approach is to remove elements without any data point from FE meshes in square domains as shown in Fig.~\ref{fig:mesh_polygon_trim}(b).

The first approach constructs a bounding polygonal domain from data directly, which adapts to random data distributions and controls the distance from the data to boundaries. However, it requires parameter tuning, including alpha parameters and mesh sizes, and may not perform well when the data is distributed in a concave polygon. In contrast, the second approach is robust and only requires the mesh size of the initial mesh, which decides the efficiency and distance from the data to boundaries. It often leads to inconsistent distance from the data to boundaries, as shown in Fig.~\ref{fig:mesh_polygon_trim}(b). It may also generate meshes with disconnected elements and the initial mesh size should be estimated by the largest distance between nearby data points.

	
Another major difference between these two approaches is the varied mesh sizes and shapes of elements. The second approach produces meshes with isosceles right triangular elements with $45^{\circ}$ or $90^{\circ}$ angles. Using the newest node bisection, the resulting refined meshes will only contain isosceles right triangles. In contrast, the first approach often leads to elements with various sizes and shapes. If angles in initial meshes are bounded between~$[\theta_{min},\theta_{max}]$, angles of resulting adaptive meshes will be bounded between $[\frac{1}{2}\theta_{min},2\theta_{max}]$ using adaptive refinement with the newest node bisection technique. While they are not as optimal as the second approach, the first approach will not lead to long and thin triangles and the deterioration of interpolation errors.


\section{Adaptive refinement and error indicators}\label{app:adaptive}

The triangular elements of FE meshes are refined using the newest node bisection method and detailed descriptions are available in~\cite{mitchell1989comparison}. A vertex of each triangle is labelled as the newest node and a triangle is divided into two children triangles by adding a new edge from the newest node to the midpoint of its opposite edge. The midpoint becomes the newest node for further refinement. The angles of the children triangles are bounded, which prevents long thin triangles and high interpolation errors.

\begin{figure} [h]
    	\centering
    	\includegraphics[scale=0.9]{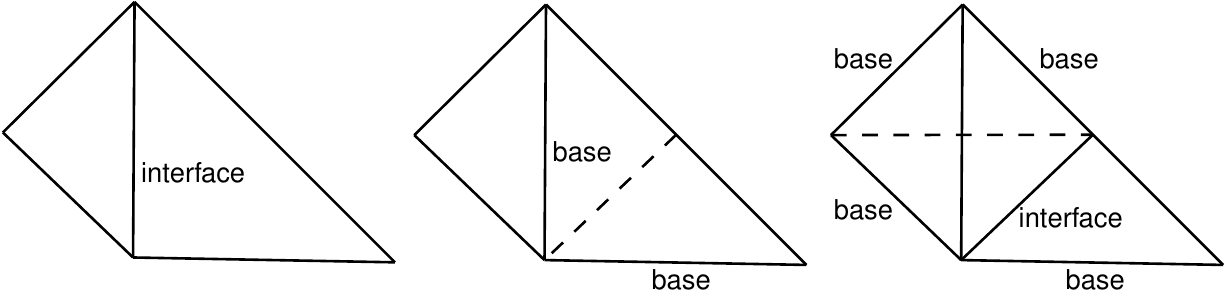}
    	\caption{Refinement of interface base-edge.}
    	\label{fig:base_interface}
    \end{figure}

The node-edge data structure is used to store the FE meshes and it is easier to focus on edges for adaptive refinement. In a triangle, the edge that sits opposite to the newest nodes is labelled as a base-edge. Any two triangles that share the same base-edge are called a triangle pair. During uniform refinement, all triangles are bisected along their base-edges and the remaining edges become base-edges for further refinement. During adaptive refinement, only a subset of the triangle pairs are bisected. Thus, the resulting FE mesh contains triangles of different sizes, which are referred to as fine and coarse triangles. The fine and coarse triangles often share common edges, which are labelled as interface base-edges. If an interface base-edge is to be bisected, the coarse triangle sharing that edge is bisected. The interface base-edge then becomes a base-edge and is bisected as illustrated in Fig.~\ref{fig:base_interface}. This prevents hanging nodes and preserves the compatibility of FE meshes.  If the coarse triangle also contains an interface base-edge, a recursive algorithm will be used to iterate through the neighbouring coarse triangles until a base-edge is found. This recursive algorithm is guaranteed to terminate.

The error indicators indicate sensitive regions of the domain and mark base-edges with high error indicator values for adaptive refinement. Four error indicators were adapted for the TPSFEM, which are auxiliary problem, residual-based, recovery-based and norm-based error indicators~\citep{fang2024error}. The auxiliary problem error indicator builds a locally more accurate approximation and approximates the error using the difference between it and $s$. The residual-based error indicator estimates the residuals of the solution directly. g The recovery-based error indicator computes improved gradients of $s$ and estimates errors by comparing it to $\nabla s$. The norm-based error indicator approximates second-order derivatives of smoother~$s$ to indicate large errors. The first two error indicators use local data and the other two do not access data directly. Out of the four error indicators, the auxiliary problem and recovery-based error indicators are compared in the numerical experiments.

The \textbf{auxiliary problem error indicator}'s local approximation is obtained by building smoother~$\hat{s}=\hat{\bm{b}}(\bm{x})^{T}\hat{\bm{c}}$ on a small subset of the domain $\Omega_e$ using a cluster of neighbouring data points~$\{(\hat{\bm{x}}_{i}, \hat{y}_{i}):i=1,2,\ldots,n_e\}$. It minimises the functional 
    \begin{equation} \label{eqn:auxi_minimiser}
		J_{\alpha}\left(\hat{\bm{c}},\hat{\bm{g}}_{1},\ldots,\hat{\bm{g}}_{d}\right) = \frac{1}{{n_e}}\sum^{{n_e}}_{i=1}\left(\hat{\bm{b}}(\hat{\bm{x}}_{i})^{T}\hat{\bm{c}}-\hat{y}_{i}\right)^{2}+\alpha \int_{{\Omega_e}} \sum_{k=1}^{d} \nabla \left(\hat{\bm{b}}(\hat{\bm{x}})^{T}\hat{\bm{g}}_{k}\right)\nabla \left(\hat{\bm{b}}(\hat{\bm{x}})^{T}\hat{\bm{g}}_{k}\right)\,d\bm{x}
    \end{equation}
subject to constraint~$\hat{L}\hat{\bm{c}} = \sum_{k=1}^{d}\hat{G}_{k}\hat{\bm{g}}_{k}$. The error norm is then calculated using the difference between~$s$ and~$\hat{s}$. The error indicator value~$\eta_{e}$ of an edge~$e$ is defined as
    \begin{equation*}
        \eta_{e}^{2}= \int_{\tau_e}\left[
        \sum_{k=1}^{d}
        \left( \frac{\partial}{\partial x_{k}}(s-\hat{s})\right)^{2}
        \right]\,d\bm{x},
    \end{equation*}
where~$\tau_e$ are triangular elements that share edge~$e$. The accuracy of auxiliary problems depends on several factors, including~$\alpha$,~$h$, size and distribution of data points and boundary condition.

The \textbf{recovery-based error indicator} estimates the error by post-processing discontinuous gradients of smoother~$s$ across inter-element boundaries. Since $s(\bm{x})=\bm{b}(\bm{x})^{T}\bm{c}$ consists of piecewise linear basis functions~$\bm{b}$, its gradient approximations~$D^{1}s$ are piecewise constants. Improve gradients~$\hat{D}^{1}s$ are obtained by~$L^{2}$-projection
    \begin{equation*} \label{eqn:recovery_improved_project}
        \int_{\Omega}b_{p}\left(\hat{D}^{1}s-D^{1}s\right)\,d\mathbf{x} = 0, \qquad p =1,\ldots,m
    \end{equation*}
	and the error is estimated as the energy norm of the difference between $D^{1}s$ and $\hat{D}^{1}s$. Its error indicator value~$\eta_{\tau_e}$ on  triangular element~$\tau_e$ is set as
    \begin{equation*}
        \eta_{\tau_e}^{2}= \int_{\tau_e}(\hat{D}^{1}s-D^{1}s)^{2}\,d\bm{x}.
    \end{equation*}	
This error indicator only uses information from the FE mesh and is computationally efficient. The descriptions of the newest node bisection and the two error indicators are in two-dimensional but they may also be extended to other dimensions.

\section{Accuracy of TPS and its derivatives}\label{app:boundary_accuracy}

A numerical experiment is conducted to evaluate the performance of the TPS for approximating a model problem and its derivatives with the two sampling strategies discussed in Sect.~\ref{sec:boundary}. The test data set is modelled by 
	\begin{equation} \label{eqn:peak}
		\begin{split}
		f(x_1,x_2) & = 3\left(1-x_1\right)^2\exp\left(-\left(x_1^2\right) - \left(x_2+1\right)^2\right) - 10\left(\frac{x_1}{5} - x_1^3 - x_2^5\right) \\
		 & \exp\left(-x_1^2-x_2^2\right)- \frac{1}{3}\exp\left(-\left(x_1+1\right)^2 - x_2^2\right)+\epsilon,
		\end{split}
	\end{equation}
where $\epsilon$ represents Gaussian noise with mean 0 and standard deviation 0.02. It consists of 10,000 data points uniformly distributed inside the $[-2.4,2.4]^2$ region. The TPS interpolant~$t$ is built with $\hat{n}$ sampled points and $\alpha$ calculated using the GCV. The quadtree sampling selects data points across the domain and boundary sampling only selects data points outside the $[-1.9,1.9]^2$ region. The accuracy of $t$ and its derivatives is measured using the RMSE, where RMSE~$ = \big(\frac{1}{n_{test}}\sum_{i=1}^{n_{test}}\big(t(\bm{x}_{i})-f(\bm{x}_{i})\big)^{2}\big)^{\frac{1}{2}}$ for $t$ and $n_{test}$ is the number of data points outside the $[-2.2,2.2]^2$ region. They are chosen to assess the approximations as they better reflect the accuracy near boundaries. The RMSE for approximations of derivatives of $f$ are defined similarly.


\begin{figure} [h]
		\centering
		\begin{subfigure}[b]{0.49\textwidth}
                \centering
			\includegraphics[scale=0.4]{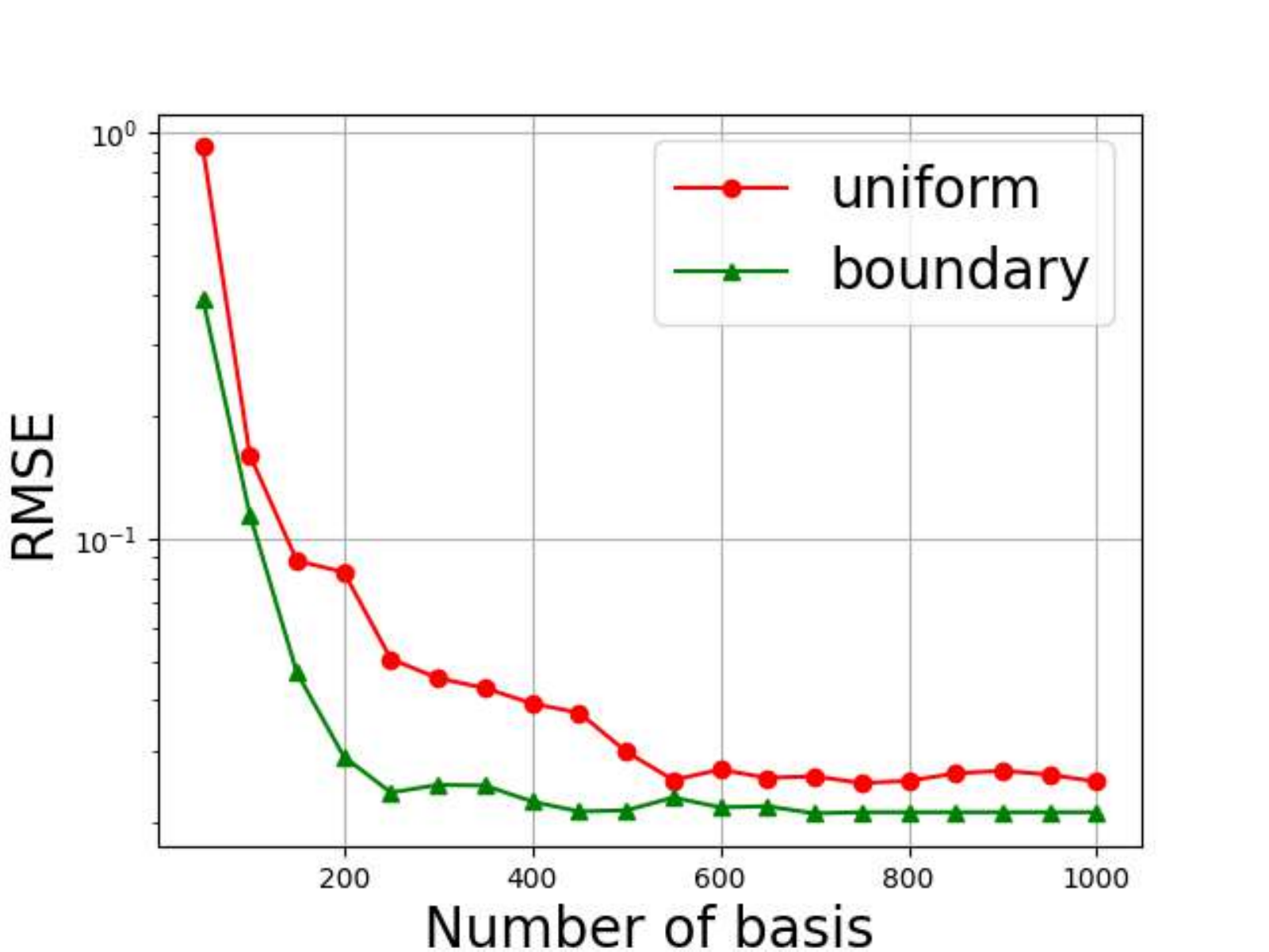}
			\caption{}
		\end{subfigure}
		\begin{subfigure}[b]{0.49\textwidth}
                \centering
			\includegraphics[scale=0.4]{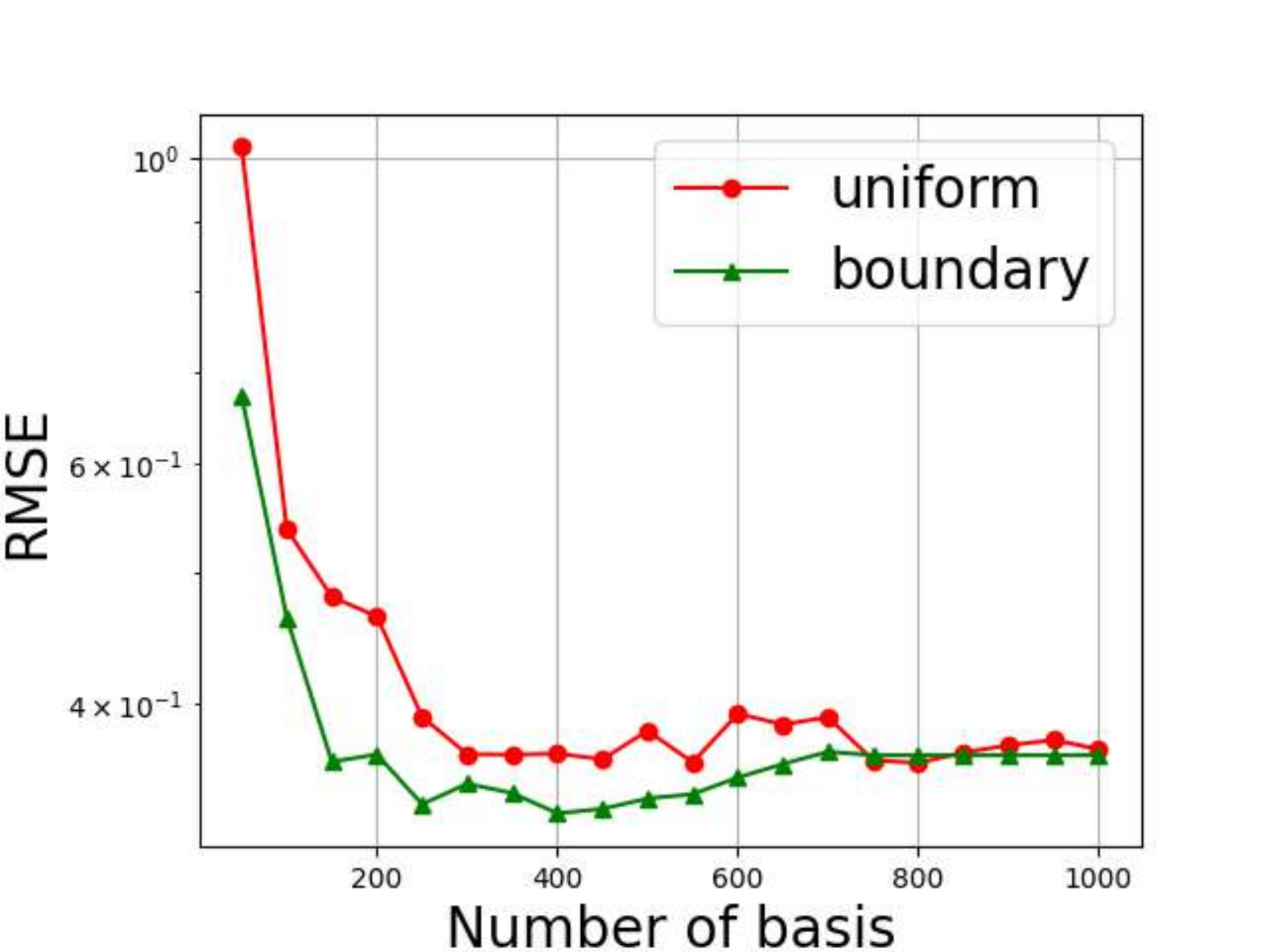}
			\caption{}
		\end{subfigure}
		\begin{subfigure}[b]{0.49\textwidth}
                \centering
			\includegraphics[scale=0.4]{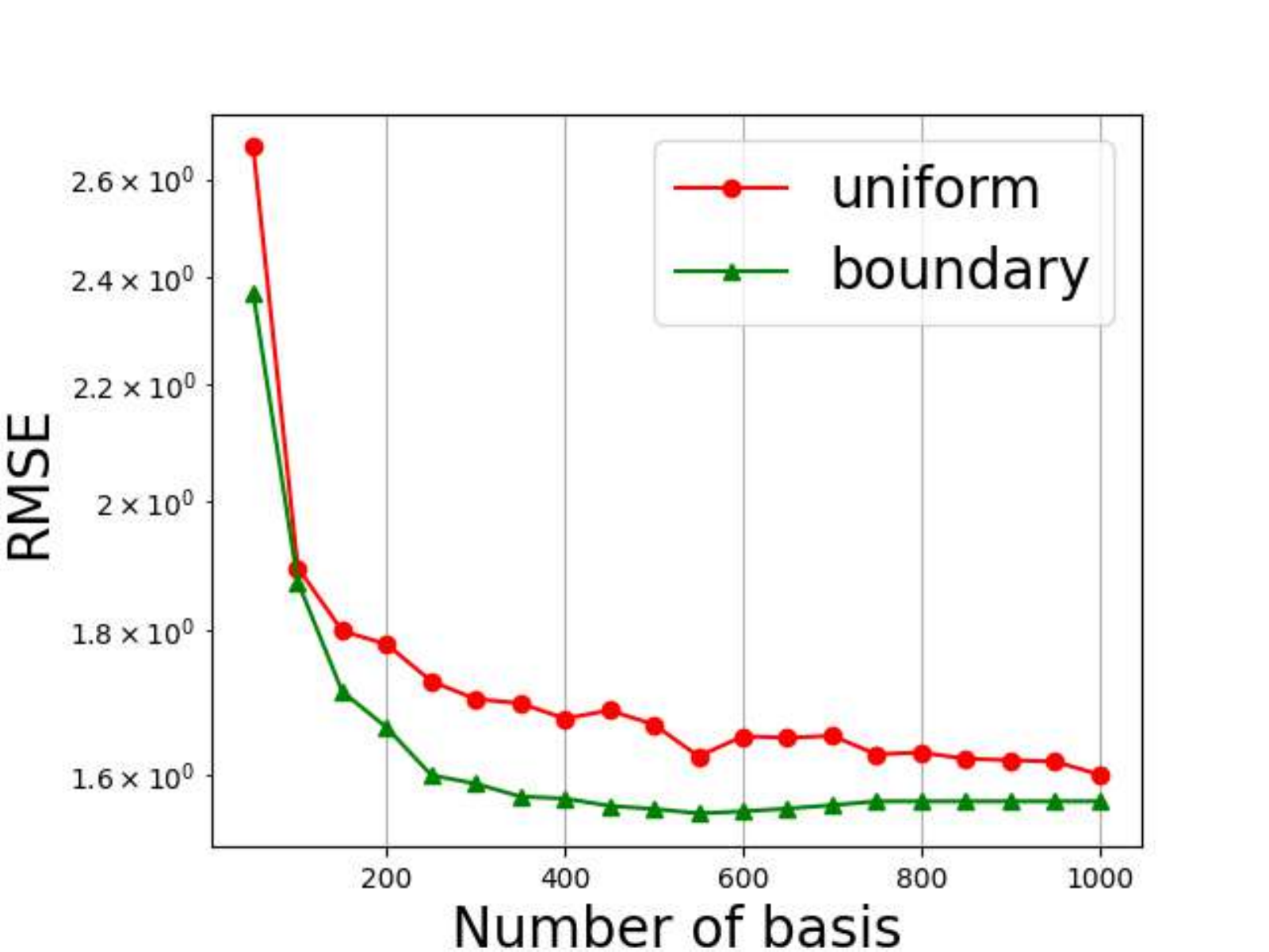}
			\caption{}
		\end{subfigure}
		\begin{subfigure}[b]{0.49\textwidth}
                \centering
			\includegraphics[scale=0.4]{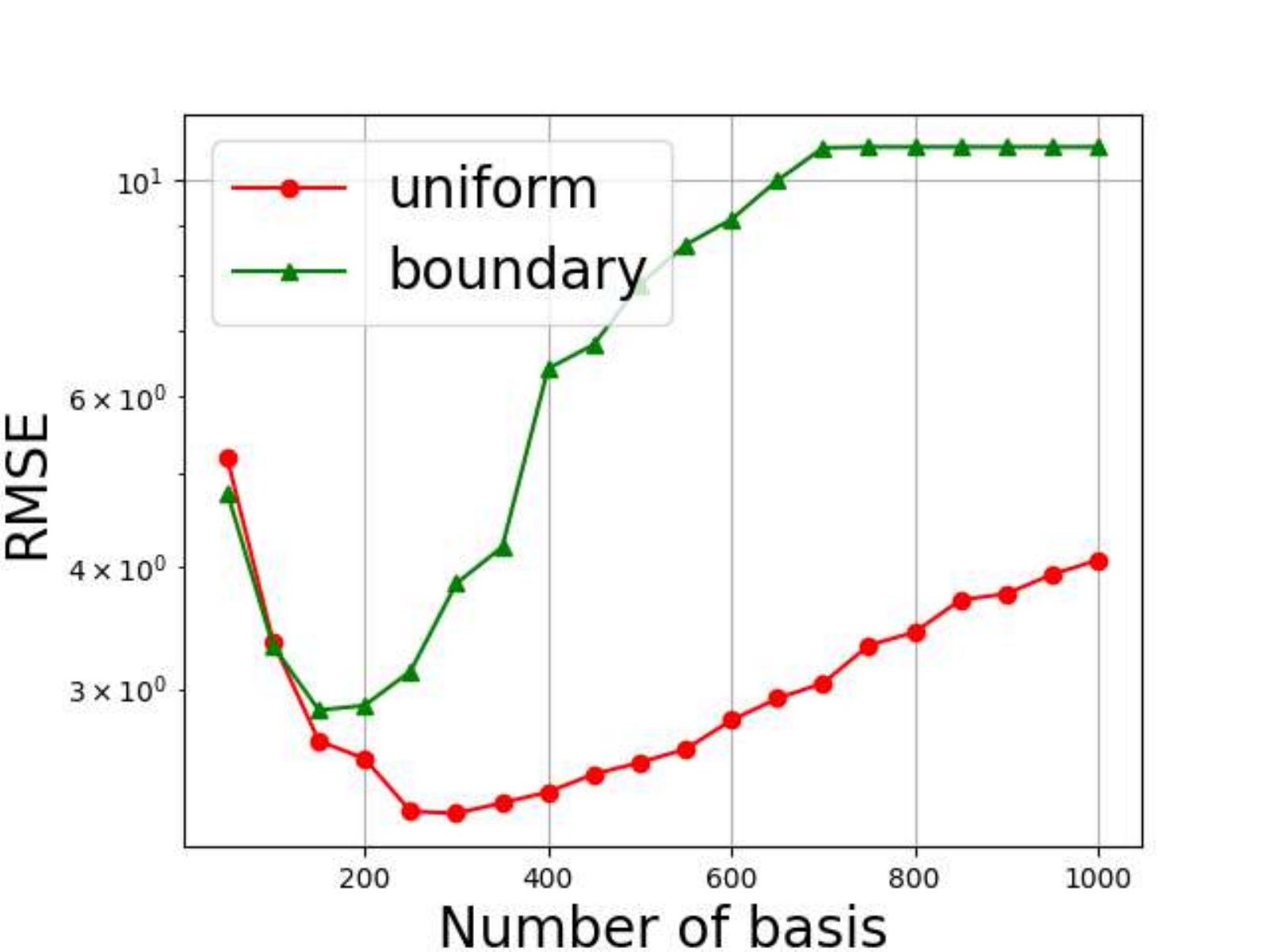}
			\caption{}
		\end{subfigure}
    \caption{RMSE of approximated (a) $f$; (b) $\frac{\partial f}{\partial x_1}$; (c) $\frac{\partial f}{\partial x_2}$; and (d) $\frac{\partial^2 f}{\partial x_1x_1}+\frac{\partial^2 f}{\partial x_2x_2}$ values using TPS kernels.}
    \label{fig:tps_errors}
\end{figure}

The convergence of the RMSE for the approximations using the two sampling strategies is shown in Fig.~\ref{fig:tps_errors}. As the number of sampled points~$\hat{n}$ increases, the RMSE for $f$, $\frac{\partial f}{\partial x_1}$, and $\frac{\partial f}{\partial x_2}$ values initially decreases and stabilises after $n_{test}$ reaches about 300 or 400 as shown in Figs.~\ref{fig:tps_errors}(a) to~\ref{fig:tps_errors}(c) for both sampling strategies. In contrast, the RMSE increases after $\hat{n}$ reaches 200 in Fig.~\ref{fig:tps_errors}(d), which shows that the approximations of second-order derivatives are sensitive to noise and less accurate than others. In addition, $\bm{w}$ may not be approximated by the second-order derivatives when they are not zero~\citep{stals2006smoothing}. The boundary sampling performs relatively better than quadtree sampling in the experiment, however, the difference in RMSE gradually diminishes as $\hat{n}$ increases.


\end{document}